\documentclass[final,3p]{elsarticle}
 \usepackage{graphics}
 \usepackage{graphicx}
 \usepackage{epsfig}
\usepackage{amssymb}
 \usepackage{amsthm}
 \usepackage{lineno}
 \usepackage{amsmath}
   \numberwithin{equation}{section}
\usepackage{mathrsfs}

\newtheorem{thm}{Theorem}[section]

\newtheorem{lem}[thm]{Lemma}

\biboptions{sort&compress,square}
\journal{some journals}
\begin{document}
\begin{frontmatter}
\author[rvt1]{Sining Wei}
\ead{weisn835@nenu.edu.cn}
\author[rvt2]{Yong Wang\corref{cor2}}
\ead{wangy581@nenu.edu.cn}

\cortext[cor2]{Corresponding author.}

\address[rvt1]{School of Data Science and Artificial Intelligence, Dongbei University of Finance and Economics, \\
Dalian, 116025, P.R.China}
\address[rvt2]{School of Mathematics and Statistics, Northeast Normal University, Changchun, 130024, P.R.China}

\title{ Conformal Perturbations of modified Novikov Operators and the Kastler-Kalau-Walze type theorem}
\begin{abstract}
In this paper, we obtain two Kastler-Kalau-Walze type theorems
for conformal perturbations of modified Novikov Operators on four-dimensional and six-dimensional compact manifolds with (respectively without)boundary.
\end{abstract}
\end{frontmatter}
\section{Introduction}
\label{1}
The noncommutative residue found in \cite{Gu,Wo} plays a prominent role in noncommutative geometry. For one-dimensional manifolds,
the noncommutative residue was discovered by Adler \cite{MA} in connection with geometric aspects of nonlinear partial differential equations. For arbitrary closed compact $n$-dimensional manifolds, the noncommutative residue was introduced by Wodzicki in \cite{Wo} using the theory of zeta functions of elliptic pseudodifferential operators. In \cite{Co1}, Connes used the noncommutative residue to derive a conformal 4-dimensional Polyakov action analogy.
Furthermore, Connes made a challenging observation that the noncommutative residue of the square of the inverse of the
Dirac operator was proportional to the Einstein-Hilbert action in \cite{Co2}.
In \cite{Ka}, Kastler gave a brute-force proof of this theorem. In \cite{KW}, Kalau and Walze proved this theorem in the
normal coordinates system simultaneously. And then, Ackermann proved that
the Wodzicki residue  of the square of the inverse of the
Dirac operator ${\rm  Wres}(D^{-2})$ in turn is essentially the second coefficient
of the heat kernel expansion of $D^{2}$ in \cite{Ac}.

In \cite{RP}, Ponge defined lower dimensional volumes of Riemannian manifolds by the Wodzicki residue. Fedosov et al. defined a noncommutative residue on Boutet de Monvel's algebra and proved that it was a unique continuous trace in \cite{FGLS}. In \cite{S}, Schrohe gave the relation between the Dixmier trace and the noncommutative residue for manifolds with boundary. In \cite{Wa3}, Wang generalized the Kastler-Kalau-Walze type theorem to the cases of 3, 4-dimensional spin manifolds with boundary and proved a Kastler-Kalau-Walze type theorem.  In \cite{Wa3,Wa4,WJ2,WJ3,WJ4}, Y. Wang and his coauthors computed the lower dimensional volumes for 5, 6, 7-dimensional spin manifolds with boundary and also got some Kastler-Kalau-Walze type theorems. In \cite{WJ5}, authors computed $\widetilde{{\rm Wres}}[(\pi^+D^{-2})\circ(\pi^+D^{-n+2})]$ for any-dimensional manifolds with boundary, and proved a general Kastler-Kalau-Walze type theorem.

In \cite{AAA}, Y. Wang proves a Kastler-Kalau-Walze type theorem for perturbations of Dirac operators on compact manifolds with or without boundary.
In \cite{WJ6}, J. Wang and Y. Wang proved two kinds of Kastler-Kalau-Walze type theorems
for conformal perturbations
of twisted Dirac operators and conformal perturbations
of signature operators by a vector bundle with a non-unitary connection on four-dimensional manifolds with  (respectively without)boundary.
In \cite{lkl}, L\'{o}pez and his collaborators introduced an elliptic differential operator which is called Novikov operator.
In \cite{SW}, we prove Kastler-Kalau-Walze-type theorems for modified Novikov operators on compact manifolds with (respectively without) a boundary.

{\bf The motivation of this paper} is to establish two Kastler-Kalau-Walze type theorems for conformal perturbations
of modified Novikov operators on manifolds with boundary. We know that the leading symbol of conformal perturbations
of modified Novikov operators is not $ic(\xi)$. This is the reason that we study the residue of conformal perturbations
of modified Novikov operators.

 This paper is organized as follows: In Section 2, we recall
some basic formulas about Boutet de Monvel's calculus and modified Novikov Operators. In Section 3, we give a Kastler-Kalau-Walze type theorem for conformal perturbations
of modified Novikov Operators on four-dimensional manifolds with boundary. In Section 4, we give a Kastler-Kalau-Walze type theorem for conformal perturbations of modified Novikov Operators on six-dimensional manifolds with boundary.
The main results are Theorem 3.4 and Theorem 4.3 in this paper.
\section{Boutet de Monvel's calculus and modified Novikov Operators}
In this section, we shall recall
some basic facts and formulas about Boutet de Monvel's calculus.
Let $$ F:L^2({\bf R}_t)\rightarrow L^2({\bf R}_v);~F(u)(v)=\int_{R} e^{-ivt}u(t)\texttt{d}t$$ denote the Fourier transformation and
$\varphi(\overline{{\bf R}^+}) =r^+\varphi({\bf R})$ (similarly define $\varphi(\overline{{\bf R}^-}$)), where $\varphi({\bf R})$
denotes the Schwartz space and
  \begin{equation}
r^{+}:C^\infty ({\bf R})\rightarrow C^\infty (\overline{{\bf R}^+});~ f\rightarrow f|\overline{{\bf R}^+};~
 \overline{{\bf R}^+}=\{x\geq0;x\in {\bf R}\}.
\end{equation}
We define $H^+=F(\varphi(\overline{{\bf R}^+}));~ H^-_0=F(\varphi(\overline{{\bf R}^-}))$ which are orthogonal to each other. We have the following
 property: $h\in H^+~$(resp.~$H^-_0)$ if and only if $h\in C^\infty({\bf R})$ which has an analytic extension to the lower (resp.~upper) complex
half-plane $\{{\rm Im}\xi<0\}~$(resp.~$\{{\rm Im}\xi>0\})$ such that for all nonnegative integer $l$,
 \begin{equation}
\frac{d^{l}h}{d\xi^l}(\xi)\sim\sum^{\infty}_{k=1}\frac{d^l}{d\xi^l}(\frac{c_k}{\xi^k})
\end{equation}
as $|\xi|\rightarrow +\infty,{\rm Im}\xi\leq0~$(resp.~${\rm Im}\xi\geq0)$.

 Let $H'$ be the space of all polynomials and $H^-=H^-_0\bigoplus H';~H=H^+\bigoplus H^-.$ Denote by $\pi^+~$(resp.~$\pi^-)$ respectively the
 projection on $H^+~$(resp.~$H^-)$. For calculations, we take $H=\widetilde H=\{$rational functions having no poles on the real axis$\}$ ($\tilde{H}$
 is a dense set in the topology of $H$). Then on $\tilde{H}$,
 \begin{equation}
\pi^+h(\xi_0)=\frac{1}{2\pi i}\lim_{u\rightarrow 0^{-}}\int_{\Gamma^+}\frac{h(\xi)}{\xi_0+iu-\xi}d\xi,
\end{equation}
where $\Gamma^+$ is a Jordan close curve included ${\rm Im}\xi>0$ surrounding all the singularities of $h$ in the upper half-plane and
$\xi_0\in {\bf R}$. Similarly, define $\pi^{'}$ on $\tilde{H}$,
 \begin{equation}
\pi'h=\frac{1}{2\pi}\int_{\Gamma^+}h(\xi)d\xi.
\end{equation}
So, $\pi'(H^-)=0$. For $h\in H\bigcap L^1({\bf R})$, $\pi'h=\frac{1}{2\pi}\int_{{\bf R}}h(v)dv$ and for $h\in H^+\bigcap L^1({\bf R})$, $\pi'h=0$.

Denote by $\mathcal{B}$ Boutet de Monvel's algebra.
For a detailed introduction to Boutet de Monvel's algebra see Boutet de
Monvel \cite{LB}, Grubb \cite{GG}, Rempel-Schulze \cite{BS} or Schrohe-Schulze \cite{EB}.
In the following we will give a review of some basic fact we need.

An operator of order $m\in {\bf Z}$ and type $d$ is a matrix
$$A=\left(\begin{array}{lcr}
  \pi^+P+G  & K  \\
   T  &  S
\end{array}\right):
\begin{array}{cc}
\   C^{\infty}(X,E_1)\\
 \   \bigoplus\\
 \   C^{\infty}(\partial{X},F_1)
\end{array}
\longrightarrow
\begin{array}{cc}
\   C^{\infty}(X,E_2)\\
\   \bigoplus\\
 \   C^{\infty}(\partial{X},F_2)
\end{array},
$$
where $X$ is a manifold with boundary $\partial X$ and
$E_1,E_2~$(resp.~$F_1,F_2)$ are vector bundles over $X~$(resp.~$\partial X
)$.~Here,~$P:C^{\infty}_0(\Omega,\overline {E_1})\rightarrow
C^{\infty}(\Omega,\overline {E_2})$ is a classical
pseudodifferential operator of order $m$ on $\Omega$, where
$\Omega$ is an open neighborhood of $X$ and
$\overline{E_i}|X=E_i~(i=1,2)$. Then $P$ has an extension:
$~{\cal{E'}}(\Omega,\overline {E_1})\rightarrow
{\cal{D'}}(\Omega,\overline {E_2})$, where
${\cal{E'}}(\Omega,\overline {E_1})$ and ${\cal{D'}}(\Omega,\overline
{E_2})$ are the dual space of $C^{\infty}(\Omega,\overline
{E_1})$ and $C^{\infty}_0(\Omega,\overline {E_2})$. Let
$e^+:C^{\infty}(X,{E_1})\rightarrow{\cal{E'}}(\Omega,\overline
{E_1})$ denote extension by zero from $X$ to $\Omega$ and
$r^+:{\cal{D'}}(\Omega,\overline{E_2})\rightarrow
{\cal{D'}}(\Omega, {E_2})$ denote the restriction from $\Omega$ to
$X$, then define
$$\pi^+P=r^+Pe^+:C^{\infty}(X,{E_1})\rightarrow {\cal{D'}}(\Omega,
{E_2}).$$
In addition, $P$ is supposed to have the
transmission property; this means that, for all $j,k,\alpha$, the
homogeneous component $p_j$ of order $j$ in the asymptotic
expansion of the
symbol $p$ of $P$ in local coordinates near the boundary satisfies:
$$\partial^k_{x_n}\partial^\alpha_{\xi'}p_j(x',0,0,+1)=
(-1)^{j-|\alpha|}\partial^k_{x_n}\partial^\alpha_{\xi'}p_j(x',0,0,-1),$$
then $\pi^+P$ maps $C^{\infty}(X,{E_1})$ into $C^{\infty}(X,{E_2})$
by Section 2.1 of \cite{Wa5}.

Let $G$, $T$ be respectively the singular Green operator
and the trace operator of order $m$ and type $d$. $K$ is a
potential operator and $S$ is a classical pseudodifferential
operator of order $m$ along the boundary (for detailed definition,
see [11]). Denote by $B^{m,d}$ the collection of all operators of
order $m$
and type $d$,  and $\mathcal{B}$ is the union over all $m$ and $d$.\\
\indent Recall $B^{m,d}$ is a Fr\'{e}chet space. The composition
of the above operator matrices yields a continuous map:
$B^{m,d}\times B^{m',d'}\rightarrow B^{m+m',{\rm max}\{
m'+d,d'\}}.$ Write $$\widetilde{A}=\left(\begin{array}{lcr}
 \pi^+P+G  & K \\
 T  &  \widetilde{S}
\end{array}\right)
\in B^{m,d},
 \widetilde{A}'=\left(\begin{array}{lcr}
\pi^+P'+G'  & K'  \\
 T'  &  \widetilde{S}'
\end{array} \right)
\in B^{m',d'}.$$ The composition $\widetilde{A}\widetilde{A}'$ is obtained by
multiplication of the matrices(for more details see [14]). For
example $\pi^+P\circ G'$ and $G\circ G'$ are singular Green
operators of type $d'$ and
$$\pi^+P\circ\pi^+P'=\pi^+(PP')+L(P,P').$$ Here $PP'$ is the usual
composition of pseudodifferential operators and $L(P,P')$ called
leftover term is a singular Green operator of type $m'+d$. For our case, $P,P'$ are classical pseudo differential operators, in other words $\pi^+P\in \mathcal{B}^{\infty}$ and $\pi^+P'\in \mathcal{B}^{\infty}$ .

In the following, write $\pi^+D^{-1}=\left(\begin{array}{lcr}
  \pi^+D^{-1}  & 0  \\
   0  &  0
\end{array}\right)$.
Let $M$ be a compact manifold with boundary $\partial M$. We assume that the metric $g^{M}$ on $M$ has
the following form near the boundary
 \begin{equation}
 g^{M}=\frac{1}{h(x_{n})}g^{\partial M}+dx _{n}^{2} ,
\end{equation}
where $g^{\partial M}$ is the metric on $\partial M$. Let $U\subset
M$ be a collar neighborhood of $\partial M$ which is diffeomorphic $\partial M\times [0,1)$. By the definition of $h(x_n)\in C^{\infty}([0,1))$
and $h(x_n)>0$, there exists $\tilde{h}\in C^{\infty}\big((-\varepsilon,1)\big)$ such that $\tilde{h}|_{[0,1)}=h$ and $\tilde{h}>0$ for some
sufficiently small $\varepsilon>0$. Then there exists a metric $\hat{g}$ on $\hat{M}=M\bigcup_{\partial M}\partial M\times
(-\varepsilon,0]$ which has the form on $U\bigcup_{\partial M}\partial M\times (-\varepsilon,0 ]$
 \begin{equation}
\hat{g}=\frac{1}{\tilde{h}(x_{n})}g^{\partial M}+dx _{n}^{2} ,
\end{equation}
such that $\hat{g}|_{M}=g$.
We fix a metric $\hat{g}$ on the $\hat{M}$ such that $\hat{g}|_{M}=g$.

Consider the $n-1$-form
$$\sigma(\xi)=\sum\limits^{n}_{j=1}(-1)^{j+1}\xi_{j}d\xi_{1}
\wedge\cdots\wedge\widehat{d\xi_{j}}\wedge\cdots\wedge d\xi_{n},$$
where the hat indicates that the corresponding factor has been omitted.

Restricted $\sigma(\xi)$ to the $n-1$-dimensional unit sphere $|\xi|=1$, $\sigma(\xi)$ gives the volume form on $|\xi|=1$.
Denoting by $|\xi'|=1$ and $\sigma(\xi')$ the $n-2$-dimensional unit sphere and the corresponding $n-2$-form.

Denote by $\mathcal{B}^{\infty}$ the algera of all operators in Boutet de Monvel's calculus (with integral order) and by $\mathcal{B}^{-\infty}$ the ideal of all smoothing operators in $\mathcal{B}^{\infty}$.
Now we recall the main theorem in \cite{FGLS}.

\begin{thm}\label{th:32}{\bf(Fedosov-Golse-Leichtnam-Schrohe)}
 Let $X$ and $\partial X$ be connected, ${\rm dim}X=n\geq3$,
 $A=\left(\begin{array}{lcr}\pi^+P+G &   K \\
T &  S    \end{array}\right)$ $\in \mathcal{B}$ , and denote by $p$, $b$ and $s$ the local symbols of $P,~G$ and $S$ respectively.
 Define:
 \begin{eqnarray}
{\rm{\widetilde{Wres}}}(A)&=&\int_X\int_{\bf |\xi|=1}{\rm{tr}}_E\left[p_{-n}(x,\xi)\right]\sigma(\xi)dx \nonumber\\
&&+2\pi\int_ {\partial X}\int_{\bf |\xi'|=1}\left\{{\rm tr}_E\left[({\rm{tr}}b_{-n})(x',\xi')\right]+{\rm{tr}}
_F\left[s_{1-n}(x',\xi')\right]\right\}\sigma(\xi')dx',
\end{eqnarray}
then

~~ a) ${\rm \widetilde{Wres}}([A,B])=0 $, for any $A,B\in\mathcal{B}$;

~~ b) It is a unique continuous trace on
$\mathcal{B}/\mathcal{B}^{-\infty}$.
\end{thm}

Secondly, we recall the definition of Novikov Operator (see details in \cite{lkl}). Let $M$ be a $n$-dimensional ($n\geq 3$) oriented compact Riemannian manifold with a Riemannian metric $g^{M}$. The de Rham derivative $d$ is a differential operator on $C^\infty(M;\wedge^*T^*M)$. Then we have the de Rham coderivative $\delta=d^*$, the symmetric operators $D=d+\delta$ and $\Delta=D^2=d\delta+\delta d$ (the Laplacian).

With more generality, we take any closed $\theta \in C^\infty(M;T^*M)$. For the sake of simplicity, we assume that $\theta$ is real. Then we have the Novikov operators defined by $\theta$, depending on $z\in \mathbb{C}$ in \cite{lkl},
\begin{eqnarray*}
d_z&=&d+z(\theta\wedge),~~\delta_z=d_z^*=\delta+\overline{z}(\theta\wedge)^*,\nonumber\\
D_z&=&d_z+\delta_z=(d+\delta)+z(\theta\wedge)+\overline{z}(\theta\wedge)^*\nonumber\\
   &=&(d+\delta)+[Rez(\theta\wedge)+Rez(\theta\wedge)^*]
                +i[Imz(\theta\wedge)-Imz(\theta\wedge)^*]\nonumber\\
   &=&(d+\delta)+Rez[\theta\wedge+(\theta\wedge)^*]
                +iImz[\theta\wedge-(\theta\wedge)^*]\nonumber\\
   &=&(d+\delta)+Rez\bar{c}(\theta)
                +iImzc(\theta),
\end{eqnarray*}
where $Rez$ is the real part of $z$, $Imz$ is the imaginary part of $z$, $\bar{c}(\theta)=(\theta)^*\wedge+(\theta\wedge)^*$,
$c(\theta)=(\theta)^*\wedge-(\theta\wedge)^*$.

For $\theta,~\theta'\in \Gamma(TM)$, we consider the modified Novikov operators. We define that
\begin{eqnarray*}
\widehat{D}_{N}&=&d+\delta+\bar{c}(\theta)+c(\theta'),~~\widehat{D}^*_{N}=d+\delta+\bar{c}(\theta)-c(\theta'),
\end{eqnarray*}
where $\bar{c}(\theta)=(\theta)^*\wedge+(\theta\wedge)^*$,
$c(\theta')=(\theta')^*\wedge-(\theta'\wedge)^*$, $\theta^*=g(\theta,\cdot)$, $(\theta')^*=g(\theta',\cdot)$.

Let $\nabla^L$ be the Levi-Civita connection about $g^M$. In the local coordinates $\{x_i; 1\leq i\leq n\}$ and the
fixed orthonormal frame $\{\widetilde{e_1},\cdots,\widetilde{e_n}\}$, the connection matrix $(\omega_{s,t})$ is defined by
\begin{equation}
\nabla^L(\widetilde{e_1},\cdots,\widetilde{e_n})= (\widetilde{e_1},\cdots,\widetilde{e_n})(\omega_{s,t}).
\end{equation}
 Let $\epsilon (\widetilde{e_j*} ),~\iota (\widetilde{e_j*} )$ be the exterior and interior multiplications respectively and $c(\widetilde{e_j})$ be the Clifford action.
Suppose that $\partial_{i}$ is a natural local frame on $TM$
and $(g^{ij})_{1\leq i,j\leq n}$ is the inverse matrix associated to the metric
matrix  $(g_{ij})_{1\leq i,j\leq n}$ on $M$. Write
\begin{equation}
c(\widetilde{e_j})=\epsilon (\widetilde{e_j*} )-\iota (\widetilde{e_j*} );~~
\bar{c}(\widetilde{e_j})=\epsilon (\widetilde{e_j*} )+\iota
(\widetilde{e_j*} ).
\end{equation}
 The modified Novikov Operators $\widehat{D}_{N}$ and $\widehat{D}^*_{N}$ are defined by
\begin{equation}
\widehat{D}_{N}=d+\delta+\bar{c}(\theta)+c(\theta')=\sum^n_{i=1}c(\widetilde{e_i})\bigg[\widetilde{e_i}+\frac{1}{4}\sum_{s,t}\omega_{s,t}
(\widetilde{e_i})\big[\bar{c}(\widetilde{e_s})\bar{c}(\widetilde{e_t})
-c(\widetilde{e_s})c(\widetilde{e_t})\big]\bigg]+\bar{c}(\theta)+c(\theta');
\end{equation}
\begin{equation}
\widehat{D}^*_{N}=d+\delta+\bar{c}(\theta)-c(\theta')=\sum^n_{i=1}c(\widetilde{e_i})\bigg[\widetilde{e_i}+\frac{1}{4}\sum_{s,t}\omega_{s,t}
(\widetilde{e_i})\big[\bar{c}(\widetilde{e_s})\bar{c}(\widetilde{e_t})
-c(\widetilde{e_s})c(\widetilde{e_t})\big]\bigg]+\bar{c}(\theta)-c(\theta').
\end{equation}

Let $g^{ij}=g(dx_{i},dx_{j})$, $\xi=\sum_{k}\limits\xi_{j}dx_{j}$ and $\nabla^L_{\partial_{i}}\partial_{j}=\sum\limits_{k}\Gamma_{ij}^{k}\partial_{k}$,  we denote that
\begin{eqnarray}
&&\sigma_{i}=-\frac{1}{4}\sum_{s,t}\omega_{s,t}
(\widetilde{e_i})c(\widetilde{e_s})c(\widetilde{e_t})
;~~~a_{i}=\frac{1}{4}\sum_{s,t}\omega_{s,t}
(\widetilde{e_i})\bar{c}(\widetilde{e_s})\bar{c}(\widetilde{e_t});\nonumber\\
&&\xi^{j}=g^{ij}\xi_{i};~~~~\Gamma^{k}=g^{ij}\Gamma_{ij}^{k};~~~~\sigma^{j}=g^{ij}\sigma_{i};
~~~~a^{j}=g^{ij}a_{i}.
\end{eqnarray}
Then the modified Novikov Operators $\widehat{D}_{N}$ and $\widehat{D}^*_{N}$ can be written as
\begin{equation}
\widehat{D}_{N}=\sum^n_{i=1}c(\widetilde{e_i})[\widetilde{e_i}+a_{i}+\sigma_{i}]
+\bar{c}(\theta)+c(\theta');
\end{equation}
\begin{equation}
\widehat{D}^*_{N}=\sum^n_{i=1}c(\widetilde{e_i})[\widetilde{e_i}+a_{i}+\sigma_{i}]
+\bar{c}(\theta)-c(\theta').
\end{equation}

By Lemma 1 in \cite{Wa4} and Lemma 2.1 in \cite{Wa3}, for any fixed point $x_0\in\partial M$, choosing the normal coordinates $U$
of $x_0$ in $\partial M$ (not in $M$). Denote by $\sigma_{l}(P)$ the $l$-order symbol of an operator $P$. By the composition formula and (2.2.11) in \cite{Wa3}, we obtain in [19, Lemma 2.6],

\begin{lem} The following identities hold:
\begin{eqnarray}
\sigma_1(\widehat{D}_{N})&=&\sigma_1(\widehat{D}^*_{N})=ic(\xi); \nonumber\\ \sigma_0(\widehat{D}_{N})&=&\frac{1}{4}\sum_{i,s,t}\omega_{s,t}(\widetilde{e_i})c(\widetilde{e_i})\bar{c}(\widetilde{e_s})\bar{c}(\widetilde{e_t})
-\frac{1}{4}\sum_{i,s,t}\omega_{s,t}(\widetilde{e_i})c(\widetilde{e_i})
c(\widetilde{e_s})c(\widetilde{e_t})+\bar{c}(\theta)+c(\theta'); \nonumber\\
\sigma_0(\widehat{D}^*_{N})&=&\frac{1}{4}\sum_{i,s,t}\omega_{s,t}(\widetilde{e_i})c(\widetilde{e_i})\bar{c}(\widetilde{e_s})\bar{c}(\widetilde{e_t})
-\frac{1}{4}\sum_{i,s,t}\omega_{s,t}(\widetilde{e_i})c(\widetilde{e_i})
c(\widetilde{e_s})c(\widetilde{e_t})+\bar{c}(\theta)-c(\theta').
\end{eqnarray}
\end{lem}

Write
 \begin{eqnarray}
D_x^{\alpha}&=&(-i)^{|\alpha|}\partial_x^{\alpha};
~\sigma(D)=p_1+p_0;
~\sigma(D^{-1})=\sum^{\infty}_{j=1}q_{-j}.
\end{eqnarray}

By the composition formula of pseudodifferential operators, we have
\begin{eqnarray}
1=\sigma(D\circ D^{-1})&=&\sum_{\alpha}\frac{1}{\alpha!}\partial^{\alpha}_{\xi}[\sigma(D)]
D^{\alpha}_{x}[\sigma(D^{-1})]\nonumber\\
&=&(p_1+p_0)(q_{-1}+q_{-2}+q_{-3}+\cdots)\nonumber\\
& &~~~+\sum_j(\partial_{\xi_j}p_1+\partial_{\xi_j}p_0)(
D_{x_j}q_{-1}+D_{x_j}q_{-2}+D_{x_j}q_{-3}+\cdots)\nonumber\\
&=&p_1q_{-1}+(p_1q_{-2}+p_0q_{-1}+\sum_j\partial_{\xi_j}p_1D_{x_j}q_{-1})+\cdots,
\end{eqnarray}
so
\begin{equation}
q_{-1}=p_1^{-1};~q_{-2}=-p_1^{-1}[p_0p_1^{-1}+\sum_j\partial_{\xi_j}p_1D_{x_j}(p_1^{-1})].
\end{equation}

By Lemma 2.2, we have some symbols of operators.
\begin{lem} The following identities hold:
\begin{eqnarray}
\sigma_{-1}(\widehat{D}^{-1}_{N})&=&\sigma_{-1}\big((\widehat{D}^*_{N})^{-1}\big)=\frac{ic(\xi)}{|\xi|^2};\nonumber\\
\sigma_{-2}(\widehat{D}^{-1}_{N})&=&\frac{c(\xi)\sigma_{0}(\widehat{D}^{-1}_{N})c(\xi)}{|\xi|^4}+\frac{c(\xi)}{|\xi|^6}\sum_jc(dx_j)
\Big[\partial_{x_j}(c(\xi))|\xi|^2-c(\xi)\partial_{x_j}(|\xi|^2)\Big] ;\nonumber\\
\sigma_{-2}\big((\widehat{D}^*_{N})^{-1}\big)&=&\frac{c(\xi)\sigma_{0}\big((\widehat{D}^*_{N})^{-1}\big)c(\xi)}{|\xi|^4}+\frac{c(\xi)}{|\xi|^6}\sum_jc(dx_j)
\Big[\partial_{x_j}(c(\xi))|\xi|^2-c(\xi)\partial_{x_j}(|\xi|^2)\Big]. \end{eqnarray}
\end{lem}

\section{ A Kastler-Kalau-Walze type theorem for four-dimensional
manifolds with boundary}

Let $M$ be $4$-dimensional compact manifolds with the boundary $\partial M$.
In the following, we will compute the more general case
$\widetilde{Wres}[\pi^{+}(f\widehat{D}^{-1}_{N}) \circ\pi^{+}\big(f^{-1}(\widehat{D}^{*}_{N})^{-1}\big)]$ for nonzero
smooth functions $f,~f^{-1}$.
An application of (3.5) and (3.6) in \cite{Wa5} shows that
\begin{eqnarray}
&&\widetilde{Wres}[\pi^{+}(f\widehat{D}^{-1}_{N}) \circ\pi^{+}\big(f^{-1}(\widehat{D}^{*}_{N})^{-1}\big)]
=Wres[f\widehat{D}^{-1}_{N}\circ f^{-1}(\widehat{D}^{*}_{N})^{-1}]+\int_{\partial M}\Phi,
\end{eqnarray}
where
 \begin{eqnarray}
\Phi&=&\int_{|\xi'|=1}\int_{-\infty}^{+\infty}\sum_{j,k=0}^{\infty}\sum \frac{(-i)^{|\alpha|+j+k+\ell}}{\alpha!(j+k+1)!}
{{\rm trace}}_{\wedge^*T^*M}\Big[\partial_{x_{n}}^{j}\partial_{\xi'}^{\alpha}\partial_{\xi_{n}}^{k}\sigma_{r}^{+}
(f\widehat{D}^{-1}_{N})(x',0,\xi',\xi_{n})\nonumber\\
&&\times\partial_{x_{n}}^{\alpha}\partial_{\xi_{n}}^{j+1}\partial_{x_{n}}^{k}\sigma_{l}
\Big(f^{-1}(\widehat{D}^{*}_{N})^{-1}\Big)(x',0,\xi',\xi_{n})\Big]
d\xi_{n}\sigma(\xi')dx' ,
\end{eqnarray}
and the sum is taken over $r-k+|\alpha|+\ell-j-1=-n=-4,r\leq-1,\ell\leq-1$.

Note that
\begin{eqnarray}
f\widehat{D}^{-1}_{N}\circ f^{-1}(\widehat{D}^{*}_{N})^{-1}=\Big(\widehat{D}^{*}_{N}\widehat{D}_{N}-\widehat{D}^{*}_{N}c(df)f^{-1}\Big)^{-1}.
\end{eqnarray}

We first establish the main theorem in this section. One has the following Lichnerowicz formula.

\begin{thm} The following equalities hold:
\begin{eqnarray}
&&\widehat{D}^*_{N}\widehat{D}_{N}-\widehat{D}^{*}_{N}c(df)f^{-1}\nonumber\\
&=&-\Big[g^{ij}(\nabla_{\partial_{i}}\nabla_{\partial_{j}}-
\nabla_{\nabla^{L}_{\partial_{i}}\partial_{j}})\Big]
-\frac{1}{8}\sum_{ijkl}R_{ijkl}\bar{c}(\widetilde{e_i})\bar{c}(\widetilde{e_j})
c(\widetilde{e_k})c(\widetilde{e_l})
+\sum_{i}c(\widetilde{e_i})\bar{c}(\nabla_{\widetilde{e_i}}^{TM}\theta)
+\frac{1}{4}s-c(\theta')\nonumber\\
&&\times\bar{c}(\theta)
+\bar{c}(\theta)c(\theta')+|\theta|^2+|\theta'|^2+\frac{1}{4}\sum_{i}
[c(\widetilde{e_{i}})c(\theta')-c(\theta')c(\widetilde{e_{i}})]^2
-g(\widetilde{e_{j}},\nabla^{TM}_{\widetilde{e_{j}}}\theta')+\sum_{i}c(e_{i})\nonumber\\
&&\times\frac{\partial_{i}(c(df)f^{-1})}{\partial x_{j}}
+\Big(\bar{c}(\theta)-c(\theta')\Big)c(df)f^{-1}
-\partial_{i}\Big(\frac{1}{2}c(\partial_{i})c(df)f^{-1}\Big)-\frac{1}{4}\sum_{i}\Big[c(e_{i})c(\theta')-c(\theta')
c(e_{i})\Big]\nonumber\\
&&\times c(e_{i})c(df)f^{-1}-\frac{1}{4}\sum_{i}c(e_{i})c(df)f^{-1}
\Big[c(e_{i})c(\theta')-c(\theta')c(e_{i})\Big],
\end{eqnarray}
where $s$ is the scalar curvature.
\end{thm}

In order to prove Theorem 3.1,
we recall the basic notions of Laplace type operators. Let $M$ be smooth compact oriented Riemannian $n$-dimensional manifolds without boundary and $V'$ be a vector bundle on $M$. Any differential operator $P$ of Laplace type has locally the form
\begin{equation}
P=-(g^{ij}\partial_i\partial_j+A^i\partial_i+B),
\end{equation}
where $\partial_{i}$  is a natural local frame on $TM$
and $(g^{ij})_{1\leq i,j\leq n}$ is the inverse matrix associated to the metric
matrix  $(g_{ij})_{1\leq i,j\leq n}$ on $M$,
 and $A^{i}$ and $B$ are smooth sections
of $\textrm{End}(V')$ on $M$ (endomorphism). If $P$ is a Laplace type
operator with the form (3.5), then there is a unique
connection $\nabla$ on $V'$ and a unique endomorphism $E$ such that
 \begin{equation}
P=-\Big[g^{ij}(\nabla_{\partial_{i}}\nabla_{\partial_{j}}-
 \nabla_{\nabla^{L}_{\partial_{i}}\partial_{j}})+E\Big],
\end{equation}
where $\nabla^{L}$ is the Levi-Civita connection on $M$. Moreover
(with local frames of $T^{*}M$ and $V'$), $\nabla_{\partial_{i}}=\partial_{i}+\omega_{i} $
and $E$ are related to $g^{ij}$, $A^{i}$ and $B$ through
 \begin{eqnarray}
&&\omega_{i}=\frac{1}{2}g_{ij}\big(A^{i}+g^{kl}\Gamma_{ kl}^{j} {\rm id}\big),\\
&&E=B-g^{ij}\big(\partial_{i}(\omega_{j})+\omega_{i}\omega_{j}-\omega_{k}\Gamma_{ ij}^{k} \big),
\end{eqnarray}
where $\Gamma_{ kl}^{j}$ is the  Christoffel coefficient of $\nabla^{L}$.

By Proposition 4.6 of \cite{wpz}, we have
\begin{equation}
(d+\delta+\bar{c}(\theta))^2=(d+\delta)^{2}
+\sum_{i}c(\widetilde{e_i})\bar{c}(\nabla_{\widetilde{e_i}}^{TM}\theta)+|\theta|^{2}.
\end{equation}
By \cite{Y}, the local expression of $(d+\delta)^{2}$ is
\begin{equation}
(d+\delta)^{2}
=-\triangle_{0}-\frac{1}{8}\sum_{ijkl}R_{ijkl}\bar{c}(\widetilde{e_i})\bar{c}(\widetilde{e_j})c(\widetilde{e_k})c(\widetilde{e_l})+\frac{1}{4}s.
\end{equation}
By \cite{Y} and \cite{Ac}, we have
\begin{equation}
-\triangle_{0}=\Delta=-g^{ij}(\nabla^L_{i}\nabla^L_{j}-\Gamma_{ij}^{k}\nabla^L_{k}).
\end{equation}
We note that
\begin{eqnarray}
\widehat{D}^*_{N}\widehat{D}_{N}&=&(d+\delta+\bar{c}(\theta))^2+(d+\delta)c(\theta')
+\bar{c}(\theta)c(\theta')-c(\theta')(d+\delta)-c(\theta')\bar{c}(\theta)
+|\theta'|^2,
\end{eqnarray}
\begin{eqnarray}
-c(\theta')(d+\delta)+(d+\delta)c(\theta')&=&\sum_{i,j}g^{i,j}\Big[c(\partial_{i})c(\theta')
-c(\theta')c(\partial_{i})\Big]\partial_{j}
-\sum_{i,j}g^{i,j}\Big[c(\theta')c(\partial_{i})\sigma_{i}+c(\theta')c(\partial_{i})a_{i}\nonumber\\
&&+c(\partial_{i})\partial_{i}(c(\theta'))+c(\partial_{i})\sigma_{i}c(\theta')
+c(\partial_{i})a_{i}c(\theta')\Big],
\end{eqnarray}
then we obtain
\begin{eqnarray}
&&\widehat{D}^*_{N}\widehat{D}_{N}-\widehat{D}^{*}_{N}c(df)f^{-1}\nonumber\\
&=&-\sum_{i,j}g^{i,j}\Big[\partial_{i}\partial_{j}
+2\sigma_{i}\partial_{j}+2a_{i}\partial_{j}-\Gamma_{i,j}^{k}\partial_{k}
+(\partial_{i}\sigma_{j})
+(\partial_{i}a_{j})
+\sigma_{i}\sigma_{j}+\sigma_{i}a_{j}+a_{i}\sigma_{j}+a_{i}a_{j} -\Gamma_{i,j}^{k}\sigma_{k}\nonumber\\
&&-\Gamma_{i,j}^{k}a_{k}\Big]
+\sum_{i,j}g^{i,j}\Big[c(\partial_{i})c(\theta')-c(\theta')c(\partial_{i})\Big]\partial_{j}
-\sum_{i,j}g^{i,j}\Big[c(\theta')c(\partial_{i})\sigma_{i}+c(\theta')c(\partial_{i})a_{i}
-c(\partial_{i})\partial_{i}(c(\theta'))\nonumber\\
&&-c(\partial_{i})\sigma_{i}c(\theta')-c(\partial_{i})a_{i}c(\theta')\Big]
-\frac{1}{8}\sum_{ijkl}R_{ijkl}\bar{c}(\widetilde{e_i})\bar{c}(\widetilde{e_j})
c(\widetilde{e_k})c(\widetilde{e_l})+\frac{1}{4}s
+\sum_{i}c(\widetilde{e_i})\bar{c}(\nabla_{\widetilde{e_i}}^{TM}\theta)+|\theta|^{2}\nonumber\\
&&+|\theta'|^2+\bar{c}(\theta)c(\theta')-c(\theta')\bar{c}(\theta)-\sum_{i,j}g^{i,j}c(\partial_{i})c(df)f^{-1}\partial_{i}
-\sum_{i,j}g^{i,j}c(\partial_{i})\frac{\partial_{i}(c(df)f^{-1})}{\partial x_{j}}-\sum_{i,j}g^{i,j}(a_{i}\nonumber\\
&&+\sigma_{i})c(df)f^{-1}-\Big(\bar{c}(\theta)-c(\theta')\Big)c(df)f^{-1}.
\end{eqnarray}
By (3.7), (3.8) and (3.14), we have
\begin{eqnarray}
(\omega_{i})_{\widehat{D}^*_{N}\widehat{D}_{N}-\widehat{D}^{*}_{N}c(df)f^{-1}}&=&\sigma_{i}+a_{i}
-\frac{1}{2}\Big[c(\partial_{i})c(\theta')-c(\theta')c(\partial_{i})\Big]+\frac{1}{2}c(\partial_{i})c(df)f^{-1}.
\end{eqnarray}
By (3.8), we have
\begin{eqnarray}
&&E_{\widehat{D}^*_{N}\widehat{D}_{N}-\widehat{D}^{*}_{N}c(df)f^{-1}}\nonumber\\
&=&-\sum_{i,j}g^{i,j}c(\partial_{i})\frac{\partial_{i}(c(df)f^{-1})}{\partial x_{j}}-\sum_{i,j}g^{i,j}(a_{i}
+\sigma_{i})c(df)f^{-1}-\Big(\bar{c}(\theta)-c(\theta')\Big)c(df)f^{-1}-
\partial_{i}\Big(\frac{1}{2}c(\partial_{i})\nonumber\\
&&\times c(df)f^{-1}\Big)-(\sigma^{i}+a^{i})\frac{1}{2}c(\partial_{i})c(df)f^{-1}+\frac{g^{ij}}{4}
\Big[c(\partial_{i})c(\theta')-c(\theta')c(\partial_{i})\Big]
c(\partial_{i})c(df)f^{-1}\nonumber\\
&&-\frac{1}{2}c(\partial_{i})c(df)f^{-1}(\sigma^{i}+a^{i})+\frac{g^{ij}}{4}c(\partial_{i})c(df)f^{-1}
\Big[c(\partial_{i})c(\theta')-c(\theta')c(\partial_{i})\Big]+\frac{1}{2}c(\partial_{i})c(df)f^{-1}\Gamma^{k}
\nonumber\\
&&
-c(\partial_{i})\sigma^{i}c(\theta')-c(\partial_{i})a^{i}c(\theta')
+\frac{1}{8}\sum_{ijkl}R_{ijkl}\bar{c}(\widetilde{e_i})\bar{c}(\widetilde{e_j})
c(\widetilde{e_k})c(\widetilde{e_l})
-\sum_{i}c(\widetilde{e_i})\bar{c}(\nabla_{\widetilde{e_i}}^{TM}\theta)
-|\theta|^{2}-|\theta'|^2\nonumber\\
&&-\frac{1}{4}s+c(\theta')\bar{c}(\theta)
+c(\theta')c(\partial_{i})\sigma^{i}+c(\theta')c(\partial_{i})a^{i}
-c(\partial_{i})\partial^{i}(c(\theta'))
+\frac{1}{2}\partial^{j}[c(\partial_{j})c(\theta')
-c(\theta')c(\partial_{j})]\nonumber\\
&&-\frac{1}{2}[c(\partial_{j})c(\theta')
-c(\theta')c(\partial_{j})](\sigma^{j}+a^{j})
-\frac{g^{ij}}{4}[c(\partial_{i})c(\theta')
-c(\theta')c(\partial_{i})][c(\partial_{j})c(\theta')
-c(\theta')c(\partial_{j})]\nonumber\\
&&-\frac{1}{2}\Gamma^{k}[c(\partial_{k})c(\theta')
-c(\theta')c(\partial_{k})]
-\bar{c}(\theta)c(\theta')
-\frac{1}{2}(\sigma^{j}+a^{j})
[c(\partial_{j})c(\theta')
-c(\theta')c(\partial_{j})].
\end{eqnarray}
For a smooth vector field $Y$ on $M$, let $c(Y)$
denote the Clifford action. Since $E$ is globally
defined on $M$, taking normal coordinates at $x_0$, we have
$\sigma^{i}(x_0)=0$, $a^{i}(x_0)=0$, $\partial^{j}[c(\partial_{j})](x_0)=0$,
$\Gamma^k(x_0)=0$, $g^{ij}(x_0)=\delta^j_i$, so that
\begin{eqnarray}
&&E_{\widehat{D}^*_{N}\widehat{D}_{N}-\widehat{D}^{*}_{N}c(df)f^{-1}}(x_0)\nonumber\\
&=&-\sum_{i}c(e_{i})\frac{\partial_{i}(c(df)f^{-1})}{\partial x_{j}}
-\Big(\bar{c}(\theta)-c(\theta')\Big)c(df)f^{-1}
+\partial_{i}\Big(\frac{1}{2}c(\partial_{i})c(df)f^{-1}\Big)+\frac{1}{4}\sum_{i}\Big[c(e_{i})c(\theta')\nonumber\\
&&-c(\theta')
c(e_{i})\Big]c(e_{i})c(df)f^{-1}+\frac{1}{4}\sum_{i}c(e_{i})c(df)f^{-1}
\Big[c(e_{i})c(\theta')-c(\theta')c(e_{i})\Big]+
\frac{1}{8}\sum_{ijkl}R_{ijkl}\bar{c}(\widetilde{e_i})\bar{c}(\widetilde{e_j})
\nonumber\\
&&\times c(\widetilde{e_k})c(\widetilde{e_l})
-\sum_{i}c(\widetilde{e_i})\bar{c}(\nabla_{\widetilde{e_i}}^{TM}\theta)
-\frac{1}{4}s+c(\theta')\bar{c}(\theta)-\bar{c}(\theta)c(\theta')
-|\theta|^2-|\theta'|^2-\frac{1}{4}\sum_{i}[c(\widetilde{e_{i}})c(\theta')\nonumber\\
&&-c(\theta')c(\widetilde{e_{i}})]^2
-\frac{1}{2}[c(\widetilde{e_{j}})\widetilde{e_{j}}(c(\theta'))
+\widetilde{e_{j}}(c(\theta'))c(\widetilde{e_{j}})]\nonumber\\
&=&-\sum_{i}c(e_{i})\frac{\partial_{i}(c(df)f^{-1})}{\partial x_{j}}
-\Big(\bar{c}(\theta)-c(\theta')\Big)c(df)f^{-1}
+\partial_{i}\Big(\frac{1}{2}c(e_{i})c(df)f^{-1}\Big)+\frac{1}{4}\sum_{i}\Big[c(e_{i})c(\theta')\nonumber\\
&&-c(\theta')
c(e_{i})\Big]c(e_{i})c(df)f^{-1}+\frac{1}{4}\sum_{i}c(e_{i})c(df)f^{-1}
\Big[c(e_{i})c(\theta')-c(\theta')c(e_{i})\Big]+\frac{1}{8}\sum_{ijkl}R_{ijkl}\bar{c}(\widetilde{e_i})\bar{c}(\widetilde{e_j})\nonumber\\
&&\times c(\widetilde{e_k})c(\widetilde{e_l})
-\sum_{i}c(\widetilde{e_i})\bar{c}(\nabla_{\widetilde{e_i}}^{TM}\theta)
-\frac{1}{4}s+c(\theta')\bar{c}(\theta)-\bar{c}(\theta)c(\theta')
-|\theta|^2-|\theta'|^2-\frac{1}{4}\sum_{i}[c(\widetilde{e_{i}})c(\theta')\nonumber\\
&&-c(\theta')c(\widetilde{e_{i}})]^2
+g(\widetilde{e_{j}}
,\nabla^{TM}_{\widetilde{e_{j}}}\theta'),
\end{eqnarray}
which, together with (3.6), yields Theorem 3.1.

The non-commutative residue of a generalized laplacian $\widetilde{\Delta}$ is expressed as by \cite{Ac}

\begin{equation}
(n-2)\Phi_{2}(\widetilde{\Delta})=(4\pi)^{-\frac{n}{2}}\Gamma(\frac{n}{2})\widetilde{res}(\widetilde{\Delta}^{-\frac{n}{2}+1}),
\end{equation}
where $\Phi_{2}(\widetilde{\Delta})$ denotes the integral over the diagonal part of the second
coefficient of the heat kernel expansion of $\widetilde{\Delta}$.
Now let $\widetilde{\Delta}=\widehat{D}^*_{N}\widehat{D}_{N}-\widehat{D}^{*}_{N}c(df)f^{-1}$. Since $\widehat{D}^*_{N}\widehat{D}_{N}-\widehat{D}^{*}_{N}c(df)f^{-1}$ is a generalized laplacian , we can suppose $\widehat{D}^*_{N}\widehat{D}_{N}-\widehat{D}^{*}_{N}c(df)f^{-1}=\Delta-E$, then we have
\begin{eqnarray}
{\rm Wres}\big(\widehat{D}^*_{N}\widehat{D}_{N}-\widehat{D}^{*}_{N}c(df)f^{-1}\big)^{-\frac{n-2}{2}}
=\frac{(n-2)(4\pi)^{\frac{n}{2}}}{(\frac{n}{2}-1)!}\int_{M}{\rm tr}(\frac{1}{6}s+E_{\widehat{D}^*_{N}\widehat{D}_{N}-\widehat{D}^{*}_{N}c(df)f^{-1}})d{\rm Vol_{M} },
\end{eqnarray}
where ${\rm Wres}$ denote the noncommutative residue.

We use ${\rm tr}$ as shorthand of ${\rm trace}$. One has the following Lemma.
\begin{lem} The following identity holds
\begin{eqnarray}
{\rm tr}\bigg[c(\partial_{i})\frac{\partial_{i}(c(df)f^{-1})}{\partial x_{i}}\bigg](x_0)
=\bigg[-f^{-1}\Delta(f)-\langle grad_{M}f,grad_{M}f^{-1}\rangle\bigg](x_0){\rm tr}[id];\nonumber\\
{\rm tr}\bigg[\partial_{i}\Big(c(\partial_{i})c(df)f^{-1}\Big)\bigg](x_0)
=\bigg[-f^{-1}\Delta(f)(x_0)-\langle grad_{M}f,grad_{M}f^{-1}\rangle\bigg](x_0){\rm tr}[id].
\end{eqnarray}
\end{lem}

Combining (3.19) with (3.20), we have
\begin{thm} For even $n$-dimensional compact oriented manifolds without boundary, the following equalities holds:
\begin{eqnarray}
&&{\rm Wres}\bigg(\widehat{D}^*_{N}\widehat{D}_{N}-\widehat{D}^{*}_{N}c(df)f^{-1}\bigg)^{-\frac{n-2}{2}}\nonumber\\
&=&\frac{(n-2)(4\pi)^{\frac{n}{2}}}{(\frac{n}{2}-1)!}\int_{M}
2^n\bigg\{\bigg(-\frac{1}{12}s-|\theta|^2+(n-2)|\theta'|^2+g(\widetilde{e_{j}},\nabla^{TM}_{\widetilde{e_{j}}}\theta')\bigg)
+\bigg[f^{-1}\Delta(f)\nonumber\\
&&+\langle grad_{M}f,grad_{M}f^{-1}\rangle\bigg]-\Big(\bar{c}(\theta)-c(\theta')\Big)c(df)f^{-1}
-\frac{1}{2}\bigg[f^{-1}\Delta(f)(x_0)\nonumber\\
&&+\langle grad_{M}f,grad_{M}f^{-1}\rangle\bigg]+\frac{1}{4}\sum_{i}\Big[c(e_{i})c(\theta')-c(\theta')
c(e_{i})\Big]c(e_{i})c(df)f^{-1}\nonumber\\
&&+\frac{1}{4}\sum_{i}c(e_{i})c(df)f^{-1}
\Big[c(e_{i})c(\theta')-c(\theta')c(e_{i})\Big]\bigg\}d{\rm Vol_{M}}.
\end{eqnarray}
where $s$ is the scalar curvature.
\end{thm}

Locally we can use Theorem 3.3 to compute the interior of $\widetilde{Wres}[\pi^{+}(f\widehat{D}^{-1}_{N}) \circ\pi^{+}\big(f^{-1}(\widehat{D}^{*}_{N})^{-1}\big)]$, we have
\begin{eqnarray}
&&\int_M\int_{|\xi|=1}{\rm
tr}_{\wedge^*T^*M}\bigg[\sigma_{-4}\big((\widehat{D}^*_{N}\widehat{D}_{N}-\widehat{D}^{*}_{N}c(df)f^{-1})^{-1}\big)\bigg]\sigma(\xi)dx\nonumber\\
&=&32\pi^2\int_{M}
\bigg\{\bigg[16g(\widetilde{e_{j}},\nabla^{TM}_{\widetilde{e_{j}}}\theta')-\frac{4}{3}s
-16|\theta|^2+32|\theta'|^2\bigg]+16\bigg[f^{-1}\Delta(f)+\langle grad_{M}f,grad_{M}f^{-1}\rangle\bigg]\nonumber\\
&&-16\Big(\bar{c}(\theta)-c(\theta')\Big)c(df)f^{-1}-8\bigg[f^{-1}\Delta(f)(x_0)+\langle grad_{M}f,grad_{M}f^{-1}\rangle\bigg]+4\sum_{i}\Big[c(e_{i})c(\theta')\nonumber\\
&&-c(\theta')
c(e_{i})\Big]c(e_{i})c(df)f^{-1}+4\sum_{i}c(e_{i})c(df)f^{-1}
\Big[c(e_{i})c(\theta')-c(\theta')c(e_{i})\Big]  \bigg\}d{\rm Vol_{M}}.
\end{eqnarray}

So we only need to compute $\int_{\partial M} \Phi$. From the remark above, now we can compute $\Phi$ (see formula (3.2) for the definition of $\Phi$). Since $n=4$, then ${\rm tr}_{\wedge^*T^*M}[{\rm id}]={\rm dim}(\wedge^*(4))=16$, since the sum is taken over $r+l-k-j-|\alpha|=-3,~~r\leq -1,l\leq-1,$ then we have the following five cases:

~\\
\noindent  {\bf case (a)~(I)}~$r=-1,~l=-1,~k=j=0,~|\alpha|=1$\\

\noindent By (3.2), we get
\begin{eqnarray}
&&{\rm case~(a)~(I)}\nonumber\\
&=&-\int_{|\xi'|=1}\int^{+\infty}_{-\infty}\sum_{|\alpha|=1}
 {\rm tr}[\partial^\alpha_{\xi'}\pi^+_{\xi_n}\sigma_{-1}(\widehat{D}^{-1}_{N})\times
 \partial^\alpha_{x'}\partial_{\xi_n}\sigma_{-1}((\widehat{D}^*_{N})^{-1})](x_0)d\xi_n\sigma(\xi')dx'-f\sum\limits_{j<n}\partial_{j}(f^{-1})\nonumber\\
 &&\times\int_{|\xi'|=1}\int^{+\infty}_{-\infty}\sum_{|\alpha|=1}{\rm tr}
\Big[\partial^{\alpha}_{\xi'}\pi^{+}_{\xi_{n}}\sigma_{-1}(\widehat{D}_{N}^{-1})
      \times\partial_{\xi_{n}}\sigma_{-1}
(\widehat{D}_{N}^{*})^{-1}\Big](x_0)d\xi_n\sigma(\xi')dx'\nonumber\\
 &=&0,
\end{eqnarray}

\noindent so {\rm case~(a)~(I)} vanishes.\\

 \noindent  {\bf case (a)~(II)}~$r=-1,~l=-1,~k=|\alpha|=0,~j=1$\\

\noindent By (3.2), we get
\begin{eqnarray}
&&{\rm case~(a)~(II)}\nonumber\\
&=&-\frac{1}{2}\int_{|\xi'|=1}\int^{+\infty}_{-\infty} {\rm
tr} [\partial_{x_n}\pi^+_{\xi_n}\sigma_{-1}(\widehat{D}^{-1}_{N})\times
\partial_{\xi_n}^2\sigma_{-1}((\widehat{D}^*_{N})^{-1})](x_0)d\xi_n\sigma(\xi')dx'\nonumber\\
      &&-\frac{1}{2}f^{-1}\partial_{x_{n}}(f)
      \int_{|\xi'|=1}\int^{+\infty}_{-\infty}{\rm tr}
\Big[\partial_{x_{n}}\sigma_{-1}(\widehat{D}_{N}^{-1})
      \times\partial^2_{\xi_{n}}\sigma_{-1}
(\widehat{D}_{N}^{*})^{-1}\Big](x_0)d\xi_n\sigma(\xi')dx'.
\end{eqnarray}
\noindent By Lemma 2.3, we have\\
\begin{eqnarray}\partial^2_{\xi_n}\sigma_{-1}\big((\widehat{D}^*_{N})^{-1}\big)(x_0)=i\left(-\frac{6\xi_nc(dx_n)+2c(\xi')}
{|\xi|^4}+\frac{8\xi_n^2c(\xi)}{|\xi|^6}\right);
\end{eqnarray}
\begin{eqnarray}
\partial_{x_n}\sigma_{-1}(\widehat{D}^{-1}_{N})(x_0)=\frac{i\partial_{x_n}c(\xi')(x_0)}{|\xi|^2}-\frac{ic(\xi)|\xi'|^2h'(0)}{|\xi|^4}.
\end{eqnarray}
By (2.3), (2.4) and the Cauchy integral formula we have
\begin{eqnarray}
\pi^+_{\xi_n}\left[\frac{c(\xi)}{|\xi|^4}\right](x_0)|_{|\xi'|=1}&=&\pi^+_{\xi_n}\left[\frac{c(\xi')+\xi_nc(dx_n)}{(1+\xi_n^2)^2}\right]\nonumber\\
&=&\frac{1}{2\pi i}{\rm lim}_{u\rightarrow
0^-}\int_{\Gamma^+}\frac{\frac{c(\xi')+\eta_nc(dx_n)}{(\eta_n+i)^2(\xi_n+iu-\eta_n)}}
{(\eta_n-i)^2}d\eta_n\nonumber\\
&=&-\frac{(i\xi_n+2)c(\xi')+ic(dx_n)}{4(\xi_n-i)^2}.
\end{eqnarray}
Similarly we have,
\begin{eqnarray}
\pi^+_{\xi_n}\left[\frac{i\partial_{x_n}c(\xi')}{|\xi|^2}\right](x_0)|_{|\xi'|=1}=\frac{\partial_{x_n}[c(\xi')](x_0)}{2(\xi_n-i)}.
\end{eqnarray}
By (3.26), then\\
\begin{eqnarray}\pi^+_{\xi_n}\partial_{x_n}\sigma_{-1}(\widehat{D}^{-1}_{N})|_{|\xi'|=1}
=\frac{\partial_{x_n}[c(\xi')](x_0)}{2(\xi_n-i)}+ih'(0)
\left[\frac{(i\xi_n+2)c(\xi')+ic(dx_n)}{4(\xi_n-i)^2}\right].
\end{eqnarray}
\noindent By the relation of the Clifford action and ${\rm tr}{AB}={\rm tr }{BA}$, we have the equalities:
\begin{eqnarray*}{\rm tr}[c(\xi')c(dx_n)]=0;~~{\rm tr}[c(dx_n)^2]=-16;~~{\rm tr}[c(\xi')^2](x_0)|_{|\xi'|=1}=-16;
\end{eqnarray*}
\begin{eqnarray}
{\rm tr}[\partial_{x_n}c(\xi')c(dx_n)]=0;~~{\rm tr}[\partial_{x_n}c(\xi')c(\xi')](x_0)|_{|\xi'|=1}=-8h'(0);~~{\rm tr}[\bar{c}(\widetilde{e_i})\bar{c}(\widetilde{e_j})c(\widetilde{e_k})c(\widetilde{e_l})]=0(i\neq j).
\end{eqnarray}
By (3.25), (3.29) and a direct computation, we have
\begin{eqnarray}
&&h'(0){\rm tr}\bigg[\frac{(i\xi_n+2)c(\xi')+ic(dx_n)}{4(\xi_n-i)^2}\times
\bigg(\frac{6\xi_nc(dx_n)+2c(\xi')}{(1+\xi_n^2)^2}
-\frac{8\xi_n^2[c(\xi')+\xi_nc(dx_n)]}{(1+\xi_n^2)^3}\bigg)
\bigg](x_0)|_{|\xi'|=1}\nonumber\\
&=&-16h'(0)\frac{-2i\xi_n^2-\xi_n+i}{(\xi_n-i)^4(\xi_n+i)^3}.
\end{eqnarray}
Similarly, we have
\begin{eqnarray}
&&-i{\rm
tr}\bigg[\bigg(\frac{\partial_{x_n}[c(\xi')](x_0)}{2(\xi_n-i)}\bigg)
\times\bigg(\frac{6\xi_nc(dx_n)+2c(\xi')}{(1+\xi_n^2)^2}-\frac{8\xi_n^2[c(\xi')+\xi_nc(dx_n)]}
{(1+\xi_n^2)^3}\bigg)\bigg](x_0)|_{|\xi'|=1}\nonumber\\
&=&-8ih'(0)\frac{3\xi_n^2-1}{(\xi_n-i)^4(\xi_n+i)^3}.
\end{eqnarray}
Then
\begin{eqnarray*}
&&-\frac{1}{2}\int_{|\xi'|=1}\int^{+\infty}_{-\infty} {\rm
tr} [\partial_{x_n}\pi^+_{\xi_n}\sigma_{-1}(\widehat{D}^{-1}_{N})\times
\partial_{\xi_n}^2\sigma_{-1}((\widehat{D}^*_{N})^{-1})](x_0)d\xi_n\sigma(\xi')dx'\nonumber\\
&=&-\int_{|\xi'|=1}\int^{+\infty}_{-\infty}\frac{4ih'(0)(\xi_n-i)^2}
{(\xi_n-i)^4(\xi_n+i)^3}d\xi_n\sigma(\xi')dx'\\
&=&-4ih'(0)\Omega_3\int_{\Gamma^+}\frac{1}{(\xi_n-i)^2(\xi_n+i)^3}d\xi_ndx'\\
&=&-4ih'(0)\Omega_32\pi i[\frac{1}{(\xi_n+i)^3}]^{'}|_{\xi_n=i}dx'\\
&=&-\frac{3}{2}\pi h'(0)\Omega_3dx',
\end{eqnarray*}
where ${\rm \Omega_{3}}$ is the canonical volume of $S^{3}.$\\

Similarly, we get
\begin{eqnarray}
&&-\frac{1}{2}f^{-1}\partial_{x_{n}}(f)
      \int_{|\xi'|=1}\int^{+\infty}_{-\infty}{\rm tr}
\Big[\partial_{x_{n}}\sigma_{-1}(\widehat{D}_{N}^{-1})
      \times\partial^2_{\xi_{n}}\sigma_{-1}
(\widehat{D}_{N}^{*})^{-1}\Big](x_0)d\xi_n\sigma(\xi')dx'\nonumber\\
&=&-2\pi i \Omega_3 f^{-1}\partial_{x_{n}}(f)dx'.
\end{eqnarray}

Then, we obtain
\begin{eqnarray}
{\rm case~(a)~(II)}=-\frac{3}{2}\pi h'(0)\Omega_3dx'-2\pi i \Omega_3 f^{-1}\partial_{x_{n}}(f)dx'.
\end{eqnarray}

\noindent  {\bf case (a)~(III)}~$r=-1,~l=-1,~j=|\alpha|=0,~k=1$\\

\noindent By (3.2), we get
\begin{eqnarray}
&&{\rm case~(a)~(III)}\nonumber\\
&=&-\frac{1}{2}\int_{|\xi'|=1}\int^{+\infty}_{-\infty}
{\rm tr} [\partial_{\xi_n}\pi^+_{\xi_n}\sigma_{-1}(\widehat{D}^{-1}_{N})\times
\partial_{\xi_n}\partial_{x_n}\sigma_{-1}\big((\widehat{D}^*_{N})^{-1}\big)](x_0)d\xi_n\sigma(\xi')dx'-\frac{1}{2}f\partial_{x_{n}}(f^{-1})\nonumber\\
&&\times      \int_{|\xi'|=1}\int^{+\infty}_{-\infty}{\rm tr}
\Big[\partial_{\xi_{n}}\pi^{+}_{\xi_{n}}\sigma_{-1}(\widehat{D}_{N}^{-1})
      \times\partial_{\xi_{n}}\sigma_{-1}\big((\widehat{D}^*_{N})^{-1}\big)\Big](x_0) d\xi_n\sigma(\xi')dx'.
\end{eqnarray}

\noindent By Lemma 2.3, we have
\begin{eqnarray}
&&\partial_{\xi_n}\partial_{x_n}\sigma_{-1}\big((\widehat{D}^*_{N})^{-1}\big)(x_0)|_{|\xi'|=1}\nonumber\\
&=&-ih'(0)\left[\frac{c(dx_n)}{|\xi|^4}-4\xi_n\frac{c(\xi')
+\xi_nc(dx_n)}{|\xi|^6}\right]-\frac{2\xi_ni\partial_{x_n}c(\xi')(x_0)}{|\xi|^4};
\end{eqnarray}
\begin{eqnarray}
\partial_{\xi_n}\pi^+_{\xi_n}\sigma_{-1}(\widehat{D}^{-1}_{N})(x_0)|_{|\xi'|=1}&=&-\frac{c(\xi')+ic(dx_n)}{2(\xi_n-i)^2}.
\end{eqnarray}
Similar to {\rm case~(a)~(II)}, we have\\
\begin{eqnarray}
&&{\rm tr}\left\{\frac{c(\xi')+ic(dx_n)}{2(\xi_n-i)^2}\times
ih'(0)\left[\frac{c(dx_n)}{|\xi|^4}-4\xi_n\frac{c(\xi')+\xi_nc(dx_n)}{|\xi|^6}\right]\right\}
=8h'(0)\frac{i-3\xi_n}{(\xi_n-i)^4(\xi_n+i)^3}
\end{eqnarray}
and
\begin{eqnarray}
{\rm tr}\left[\frac{c(\xi')+ic(dx_n)}{2(\xi_n-i)^2}\times
\frac{2\xi_ni\partial_{x_n}c(\xi')(x_0)}{|\xi|^4}\right]
=\frac{-8ih'(0)\xi_n}{(\xi_n-i)^4(\xi_n+i)^2}.
\end{eqnarray}
So we have
\begin{eqnarray}
&&-\frac{1}{2}\int_{|\xi'|=1}\int^{+\infty}_{-\infty}
{\rm tr} [\partial_{\xi_n}\pi^+_{\xi_n}\sigma_{-1}(\widehat{D}^{-1}_{N})\times
\partial_{\xi_n}\partial_{x_n}\sigma_{-1}\big((\widehat{D}^*_{N})^{-1}\big)](x_0)d\xi_n\sigma(\xi')dx'\nonumber\\
&=&-\int_{|\xi'|=1}\int^{+\infty}_{-\infty}\frac{h'(0)4(i-3\xi_n)}
{(\xi_n-i)^4(\xi_n+i)^3}d\xi_n\sigma(\xi')dx'
-\int_{|\xi'|=1}\int^{+\infty}_{-\infty}\frac{h'(0)4i\xi_n}
{(\xi_n-i)^4(\xi_n+i)^2}d\xi_n\sigma(\xi')dx'\nonumber\\
&=&-h'(0)\Omega_3\frac{2\pi i}{3!}[\frac{4(i-3\xi_n)}{(\xi_n+i)^3}]^{(3)}|_{\xi_n=i}dx'+h'(0)\Omega_3\frac{2\pi i}{3!}[\frac{4i\xi_n}{(\xi_n+i)^2}]^{(3)}|_{\xi_n=i}dx'\nonumber\\
&=&\frac{3}{2}\pi h'(0)\Omega_3dx'.
\end{eqnarray}

Similarly, we get
\begin{eqnarray}
&&-\frac{1}{2}f\partial_{x_{n}}(f^{-1})\int_{|\xi'|=1}\int^{+\infty}_{-\infty}{\rm tr}
\Big[\partial_{\xi_{n}}\pi^{+}_{\xi_{n}}\sigma_{-1}(\widehat{D}_{N}^{-1})
      \times\partial_{\xi_{n}}\sigma_{-1}\big((\widehat{D}^*_{N})^{-1}\big)\Big](x_0) d\xi_n\sigma(\xi')dx'\nonumber\\
&=&2\pi i \Omega_3 f\partial_{x_{n}}(f^{-1})dx'.
\end{eqnarray}

Then, we obtain
\begin{eqnarray}
{\rm case~(a)~(III)}=\frac{3}{2}\pi h'(0)\Omega_3dx'+2\pi i \Omega_3 f\partial_{x_{n}}(f^{-1})dx'.
\end{eqnarray}

\noindent  {\bf case (b)}~$r=-2,~l=-1,~k=j=|\alpha|=0$\\

\noindent By (3.2), we get
\begin{eqnarray}
{\rm case~(b)}
&=&-i\int_{|\xi'|=1}\int^{+\infty}_{-\infty}{\rm tr} [\pi^+_{\xi_n}\sigma_{-2}(f\widehat{D}^{-1}_{N})\times
\partial_{\xi_n}\sigma_{-1}(f^{-1}(\widehat{D}^*_{N})^{-1})](x_0)d\xi_n\sigma(\xi')dx'\nonumber\\
&=&-i\int_{|\xi'|=1}\int^{+\infty}_{-\infty}{\rm tr} [\pi^+_{\xi_n}\sigma_{-2}(\widehat{D}^{-1}_{N})\times
\partial_{\xi_n}\sigma_{-1}\big((\widehat{D}^*_{N})^{-1}\big)](x_0)d\xi_n\sigma(\xi')dx'.
\end{eqnarray}
 By Lemma 2.3, we have
\begin{eqnarray}
\sigma_{-2}(\widehat{D}^{-1}_{N})(x_0)=\frac{c(\xi)\sigma_{0}(\widehat{D}_{N})(x_0)c(\xi)}{|\xi|^4}+\frac{c(\xi)}{|\xi|^6}c(dx_n)
\bigg[\partial_{x_n}[c(\xi')](x_0)|\xi|^2-c(\xi)h'(0)|\xi|^2_{\partial
M}\bigg],
\end{eqnarray}
where
\begin{eqnarray}
\sigma_{0}(\widehat{D}_{N})&=&\frac{1}{4}\sum_{s,t,i}\omega_{s,t}(\widetilde{e_i})
c(\widetilde{e_i})\bar{c}(\widetilde{e_s})\bar{c}(\widetilde{e_t})
-\frac{1}{4}\sum_{s,t,i}\omega_{s,t}(\widetilde{e_i})
c(\widetilde{e_i})c(\widetilde{e_s})c(\widetilde{e_t}))+\bar{c}(\theta)+c(\theta').
\end{eqnarray}
We denote
$b_0^{1}(x_0)=\frac{1}{4}\sum_{s,t,i}\omega_{s,t}(\widetilde{e_i})(x_{0})c(\widetilde{e_i})\bar{c}(\widetilde{e_s})\bar{c}(\widetilde{e_t});~
b_0^{2}(x_0)=-\frac{1}{4}\sum_{s,t,i}\omega_{s,t}(\widetilde{e_i})(x_{0})c(\widetilde{e_i})c(\widetilde{e_s})c(\widetilde{e_t})).$

Then
\begin{eqnarray}
&&\pi^+_{\xi_n}\sigma_{-2}(\widehat{D}^{-1}_{N}(x_0))|_{|\xi'|=1}\nonumber\\
&=&\pi^+_{\xi_n}\Big[\frac{c(\xi)b_0^{1}(x_0)c(\xi)}{(1+\xi_n^2)^2}\Big]+\pi^+_{\xi_n}
\Big[\frac{c(\xi)(\bar{c}(\theta)+c(\theta'))(x_0)c(\xi)}{(1+\xi_n^2)^2}\Big]
\nonumber\\
&&+\pi^+_{\xi_n}\Big[\frac{c(\xi)b_0^{2}(x_0)c(\xi)+c(\xi)c(dx_n)\partial_{x_n}[c(\xi')](x_0)}{(1+\xi_n^2)^2}-h'(0)\frac{c(\xi)c(dx_n)c(\xi)}{(1+\xi_n^{2})^3}\Big].
\end{eqnarray}
By direct calculation, we have
\begin{eqnarray}
\pi^+_{\xi_n}\Big[\frac{c(\xi)b_0^{1}(x_0)c(\xi)}{(1+\xi_n^2)^2}\Big]&=&\pi^+_{\xi_n}\Big[\frac{c(\xi')b_0^{1}(x_0)c(\xi')}{(1+\xi_n^2)^2}\Big]
+\pi^+_{\xi_n}\Big[ \frac{\xi_nc(\xi')b_0^{1}(x_0)c(dx_{n})}{(1+\xi_n^2)^2}\Big]\nonumber\\
&&+\pi^+_{\xi_n}\Big[\frac{\xi_nc(dx_{n})b_0^{1}(x_0)c(\xi')}{(1+\xi_n^2)^2}\Big]
+\pi^+_{\xi_n}\Big[\frac{\xi_n^{2}c(dx_{n})b_0^{1}(x_0)c(dx_{n})}{(1+\xi_n^2)^2}\Big]\nonumber\\
&=&-\frac{c(\xi')b_0^{1}(x_0)c(\xi')(2+i\xi_{n})}{4(\xi_{n}-i)^{2}}
+\frac{ic(\xi')b_0^{1}(x_0)c(dx_{n})}{4(\xi_{n}-i)^{2}}\nonumber\\
&&+\frac{ic(dx_{n})b_0^{1}(x_0)c(\xi')}{4(\xi_{n}-i)^{2}}
+\frac{-i\xi_{n}c(dx_{n})b_0^{1}(x_0)c(dx_{n})}{4(\xi_{n}-i)^{2}}.
\end{eqnarray}
Since
\begin{eqnarray}
c(dx_n)b_0^{1}(x_0)
&=&-\frac{1}{4}h'(0)\sum^{n-1}_{i=1}c(\widetilde{e_i})
\bar{c}(\widetilde{e_i})c(\widetilde{e_n})\bar{c}(\widetilde{e_n}),
\end{eqnarray}
then by the relation of the Clifford action and ${\rm tr}{AB}={\rm tr }{BA}$,  we have the equalities:
\begin{eqnarray}
{\rm tr}[c(\widetilde{e_i})
\bar{c}(\widetilde{e_i})c(\widetilde{e_n})
\bar{c}(\widetilde{e_n})]&=&0~~(i<n);~~
{\rm tr}[b_0^{1}c(dx_n)]=0;~~{\rm tr }[\bar{c}(\theta)c(dx_{n})]=0;\nonumber\\
~~{\rm tr }[c(\theta')c(dx_{n})]&=&-16g(\theta',dx_n);
~~{\rm tr}[\bar {c}(\xi')\bar {c}(dx_n)]=0.
\end{eqnarray}
Since
\begin{eqnarray}
\partial_{\xi_n}\sigma_{-1}((\widehat{D}^*_{N})^{-1})=\partial_{\xi_n}q_{-1}(x_0)|_{|\xi'|=1}=i\left[\frac{c(dx_n)}{1+\xi_n^2}-\frac{2\xi_nc(\xi')+2\xi_n^2c(dx_n)}{(1+\xi_n^2)^2}\right],
\end{eqnarray}
By (3.47) and (3.50), we have
\begin{eqnarray}
&&{\rm tr }[\pi^+_{\xi_n}\Big[\frac{c(\xi)b_0^{1}(x_0)c(\xi)}{(1+\xi_n^2)^2}\Big]
\times\partial_{\xi_n}\sigma_{-1}((\widehat{D}^*_{N})^{-1})(x_0)]|_{|\xi'|=1}\nonumber\\
&=&\frac{1}{2(1+\xi_n^2)^2}{\rm tr }[c(\xi')b_0^{1}(x_0)]
+\frac{i}{2(1+\xi_n^2)^2}{\rm tr }[c(dx_n)b_0^{1}(x_0)]\nonumber\\
&=&\frac{1}{2(1+\xi_n^2)^2}{\rm tr }[c(\xi')b_0^{1}(x_0)].
\end{eqnarray}
We note that $i<n,~\int_{|\xi'|=1}\{\xi_{i_{1}}\xi_{i_{2}}\cdots\xi_{i_{2d+1}}\}\sigma(\xi')=0$,
so ${\rm tr }[c(\xi')b_0^{1}(x_0)]$ has no contribution for computing {\rm case~(b)}.

By direct calculation we have
\begin{eqnarray}
\pi^+_{\xi_n}\Big[\frac{c(\xi)b_0^{2}(x_0)c(\xi)+c(\xi)c(dx_n)\partial_{x_n}[c(\xi')](x_0)}{(1+\xi_n^2)^2}\Big]-h'(0)\pi^+_{\xi_n}\Big[\frac{c(\xi)c(dx_n)c(\xi)}{(1+\xi_n)^3}\Big]:= B_1-B_2,
\end{eqnarray}
where
\begin{eqnarray}
B_1&=&\frac{-1}{4(\xi_n-i)^2}[(2+i\xi_n)c(\xi')b_0^{2}(x_0)c(\xi')+i\xi_nc(dx_n)b_0^{2}(x_0)c(dx_n)\nonumber\\
&&+(2+i\xi_n)c(\xi')c(dx_n)\partial_{x_n}c(\xi')+ic(dx_n)b_0^{2}(x_0)c(\xi')
+ic(\xi')b_0^{2}(x_0)c(dx_n)-i\partial_{x_n}c(\xi')]
\end{eqnarray}
and
\begin{eqnarray}
B_2&=&\frac{h'(0)}{2}\left[\frac{c(dx_n)}{4i(\xi_n-i)}+\frac{c(dx_n)-ic(\xi')}{8(\xi_n-i)^2}
+\frac{3\xi_n-7i}{8(\xi_n-i)^3}[ic(\xi')-c(dx_n)]\right].
\end{eqnarray}
By (3.50) and (3.54), we have
\begin{eqnarray}{\rm tr }[B_2\times\partial_{\xi_n}\sigma_{-1}((\widehat{D}^*_{N})^{-1})]|_{|\xi'|=1}
&=&\frac{i}{2}h'(0)\frac{-i\xi_n^2-\xi_n+4i}{4(\xi_n-i)^3(\xi_n+i)^2}{\rm tr}[{\rm id}]\nonumber\\
&=&8ih'(0)\frac{-i\xi_n^2-\xi_n+4i}{4(\xi_n-i)^3(\xi_n+i)^2}.
\end{eqnarray}
By (3.50) and (3.53), we have
\begin{eqnarray}{\rm tr }[B_1\times\partial_{\xi_n}\sigma_{-1}((\widehat{D}^*_{N})^{-1})]|_{|\xi'|=1}=
\frac{-8ic_0}{(1+\xi_n^2)^2}+2h'(0)\frac{\xi_n^2-i\xi_n-2}{(\xi_n-i)(1+\xi_n^2)^2},
\end{eqnarray}
where $b_0^{2}=c_0c(dx_n)$ and $c_0=-\frac{3}{4}h'(0)$.\\

By (3.55) and (3.56), we have
\begin{eqnarray}
&&-i\int_{|\xi'|=1}\int^{+\infty}_{-\infty}{\rm tr} [(B_1-B_2)\times
\partial_{\xi_n}\sigma_{-1}\bigg((\widehat{D}^*_{N})^{-1}\bigg)](x_0)d\xi_n\sigma(\xi')dx'\nonumber\\
&=&-\Omega_3\int_{\Gamma^+}\frac{8c_0(\xi_n-i)+ih'(0)}{(\xi_n-i)^3(\xi_n+i)^2}d\xi_ndx'\nonumber\\
&=&\frac{9}{2}\pi h'(0)\Omega_3dx'.
\end{eqnarray}
Similar to (3.51), we have
\begin{eqnarray}
&&{\rm tr }[\pi^+_{\xi_n}\Big[\frac{c(\xi)\bar{c}(\theta)(x_0)c(\xi)}{(1+\xi_n^2)^2}\Big]
\times\partial_{\xi_n}\sigma_{-1}((\widehat{D}^*_{N})^{-1})(x_0)]|_{|\xi'|=1}\nonumber\\
&=&\frac{1}{2(1+\xi_n^2)^2}{\rm tr }[c(\xi')\bar{c}(\theta)(x_0)]
+\frac{i}{2(1+\xi_n^2)^2}{\rm tr }[c(dx_n)\bar{c}(\theta)(x_0)]\nonumber\\
&=&\frac{1}{2(1+\xi_n^2)^2}{\rm tr }[c(\xi')\bar{c}(\theta)(x_0)].
\end{eqnarray}
Similar to (3.51), we have
\begin{eqnarray}
{\rm tr }[\pi^+_{\xi_n}\Big[\frac{c(\xi)c(\theta')(x_0)c(\xi)}{(1+\xi_n^2)^2}\Big]
\times\partial_{\xi_n}\sigma_{-1}((\widehat{D}^*_{N})^{-1})(x_0)]|_{|\xi'|=1}
&=&\frac{i}{2(1+\xi_n^2)^2}{\rm tr }[c(dx_n)c(\theta')(x_0)].
\end{eqnarray}
By (3.58) and (3.59), we have
\begin{eqnarray}
&&-i\int_{|\xi'|=1}\int^{+\infty}_{-\infty}{\rm tr} \bigg[\pi^+_{\xi_n}
\Big[\frac{c(\xi)(\bar{c}(\theta)+c(\theta'))c(\xi)}{(1+\xi_n^2)^2}\Big]\times
\partial_{\xi_n}\sigma_{-1}((\widehat{D}^*_{N})^{-1})\bigg](x_0)d\xi_n\sigma(\xi')dx'\nonumber\\
&=&\frac{\pi}{4}{\rm tr}[c(dx_n)c(\theta')]\Omega_3dx'\nonumber\\
&=&-4\pi g(\theta',dx_n)\Omega_3dx'.
\end{eqnarray}
By (3.51), (3.57) and (3.60), we have
\begin{eqnarray}
{\rm case~(b)}=\frac{9}{2}\pi h'(0)\Omega_3dx'-4\pi g(\theta',dx_n)\Omega_3dx'.
\end{eqnarray}

\noindent {\bf  case (c)}~$r=-1,~l=-2,~k=j=|\alpha|=0$\\
By (3.2), we get
\begin{eqnarray}
{\rm case~ (c)}&=&-i\int_{|\xi'|=1}\int^{+\infty}_{-\infty}{\rm tr} [\pi^+_{\xi_n}\sigma_{-1}(f\widehat{D}^{-1}_{N})\times
\partial_{\xi_n}\sigma_{-2}(f^{-1}(\widehat{D}^*_{N})^{-1})](x_0)d\xi_n\sigma(\xi')dx'\nonumber\\
&=&-i\int_{|\xi'|=1}\int^{+\infty}_{-\infty}{\rm tr} [\pi^+_{\xi_n}\sigma_{-1}(\widehat{D}^{-1}_{N})\times
\partial_{\xi_n}\sigma_{-2}\big((\widehat{D}^*_{N})^{-1}\big)](x_0)d\xi_n\sigma(\xi')dx'.
\end{eqnarray}
By (2.3) and (2.4), Lemma 2.3, we have
\begin{equation}
\pi^+_{\xi_n}\sigma_{-1}(\widehat{D}^{-1}_{N})|_{|\xi'|=1}=\frac{c(\xi')+ic(dx_n)}{2(\xi_n-i)}.
\end{equation}

Since
\begin{equation}
\sigma_{-2}((\widehat{D}^*_{N})^{-1})(x_0)=\frac{c(\xi)\sigma_{0}(\widehat{D}^*_{N})(x_0)c(\xi)}{|\xi|^4}+\frac{c(\xi)}{|\xi|^6}c(dx_n)
\bigg[\partial_{x_n}[c(\xi')](x_0)|\xi|^2-c(\xi)h'(0)|\xi|^2_{\partial
M}\bigg],
\end{equation}
where
\begin{eqnarray}
\sigma_{0}(\widehat{D}^*_{N})(x_0)&=&\frac{1}{4}\sum_{s,t,i}\omega_{s,t}(\widetilde{e_i})(x_{0})c(\widetilde{e_i})\bar{c}(\widetilde{e_s})\bar{c}(\widetilde{e_t})
-\frac{1}{4}\sum_{s,t,i}\omega_{s,t}(\widetilde{e_i})(x_{0})c(\widetilde{e_i})
c(\widetilde{e_s})c(\widetilde{e_t}))+\big(\bar{c}(\theta)-c(\theta')\big)(x_{0})\nonumber\\
&=&b_0^{1}(x_0)+b_0^{2}(x_0)+\big(\bar{c}(\theta)-c(\theta')\big)(x_{0}),
\end{eqnarray}
then
\begin{eqnarray}
&&\partial_{\xi_n}\sigma_{-2}((\widehat{D}^*_{N})^{-1})(x_0)|_{|\xi'|=1}\nonumber\\
&=&
\partial_{\xi_n}\bigg\{\frac{c(\xi)[b_0^{1}(x_0)+b_0^{2}(x_0)
+\big(\bar{c}(\theta)-c(\theta')\big)(x_{0})]c(\xi)}{|\xi|^4}+\frac{c(\xi)}{|\xi|^6}c(dx_n)\nonumber\\
&&\times\big[\partial_{x_n}[c(\xi')](x_0)|\xi|^2-c(\xi)h'(0)\big]\bigg\}\nonumber\\
&=&\partial_{\xi_n}\bigg\{\frac{c(\xi)b_0^{1}(x_0)]c(\xi)}{|\xi|^4}\bigg\}+\partial_{\xi_n}\bigg\{\frac{c(\xi)}{|\xi|^6}
c(dx_n)[\partial_{x_n}[c(\xi')](x_0)|\xi|^2-c(\xi)h'(0)]\bigg\}\nonumber\\
&&+\partial_{\xi_n}\bigg\{\frac{c(\xi)b_0^{2}(x_0)c(\xi)}{|\xi|^4}\bigg\}
+\partial_{\xi_n}\bigg\{\frac{c(\xi)\big(\bar{c}(\theta)-c(\theta')\big)(x_{0})c(\xi)}{|\xi|^4}\bigg\}.
\end{eqnarray}
By direct calculation, we have
\begin{eqnarray}
\partial_{\xi_n}\bigg\{\frac{c(\xi)b_0^{1}(x_0)c(\xi)}{|\xi|^4}\bigg\}=\frac{c(dx_n)b_0^{1}(x_0)c(\xi)}{|\xi|^4}
+\frac{c(\xi)b_0^{1}(x_0)c(dx_n)}{|\xi|^4}
-\frac{4\xi_n c(\xi)b_0^{1}(x_0)c(\xi)}{|\xi|^6};
\end{eqnarray}
\begin{eqnarray}
\partial_{\xi_n}\bigg\{\frac{c(\xi)(\bar{c}(\theta)-c(\theta'))(x_{0})c(\xi)}{|\xi|^4}\bigg\}
&=&\frac{c(dx_n)(\bar{c}(\theta)-c(\theta'))(x_{0})c(\xi)}{|\xi|^4}
+\frac{c(\xi)(\bar{c}(\theta)-c(\theta'))(x_{0})c(dx_n)}{|\xi|^4}\nonumber\\
&&-\frac{4\xi_n c(\xi)(\bar{c}(\theta)-c(\theta'))(x_{0})c(\xi)}{|\xi|^4},
\end{eqnarray}
and
\begin{eqnarray}
&&\partial_{\xi_n}\bigg\{\frac{c(\xi)}{|\xi|^6}
c(dx_n)\big[\partial_{x_n}[c(\xi')](x_0)|\xi|^2-c(\xi)h'(0)\big]\bigg\}+\partial_{\xi_n}\bigg\{\frac{c(\xi)b_0^{2}(x_0)c(\xi)}{|\xi|^4}\bigg\}\nonumber\\
&=&\frac{1}{(1+\xi_n^2)^3}\bigg[(2\xi_n-2\xi_n^3)c(dx_n)b_0^{2}c(dx_n)
+(1-3\xi_n^2)c(dx_n)b_0^{2}c(\xi')+(1-3\xi_n^2)c(\xi')b_0^{2}c(dx_n)\nonumber\\
&&
-4\xi_nc(\xi')b_0^{2}c(\xi')
+(3\xi_n^2-1)\partial_{x_n}c(\xi')-4\xi_nc(\xi')c(dx_n)\partial_{x_n}c(\xi')+2h'(0)c(\xi')+2h'(0)\xi_nc(dx_n)\bigg]\nonumber\\
&&+6\xi_nh'(0)\frac{c(\xi)c(dx_n)c(\xi)}{(1+\xi^2_n)^4}.
\end{eqnarray}
By (3.63) and (3.67), we have
\begin{eqnarray}
&&{\rm tr}\bigg[\pi^+_{\xi_n}\sigma_{-1}(\widehat{D}^{-1}_{N})\times
\partial_{\xi_n}\frac{c(\xi)b_0^{1}c(\xi)}
{|\xi|^4}\bigg](x_0)|_{|\xi'|=1}\nonumber\\
&=&\frac{-1}{(\xi-i)(\xi+i)^3}{\rm tr}[c(\xi')b_0^{1}(x_0)]
+\frac{i}{(\xi-i)(\xi+i)^3}{\rm tr}[c(dx_n)b_0^{1}(x_0)].
\end{eqnarray}
By (3.48) and (3.49), we have
\begin{eqnarray}
{\rm tr}\bigg[\pi^+_{\xi_n}\sigma_{-1}(\widehat{D}^{-1}_{N})\times
\partial_{\xi_n}\frac{c(\xi)b_0^{1}c(\xi)}
{|\xi|^4}\bigg](x_0)|_{|\xi'|=1}
=\frac{-1}{(\xi-i)(\xi+i)^3}{\rm tr}[c(\xi')b_0^{1}(x_0)].
\end{eqnarray}
We note that $i<n,~\int_{|\xi'|=1}\{\xi_{i_{1}}\xi_{i_{2}}\cdots\xi_{i_{2d+1}}\}\sigma(\xi')=0$,
so ${\rm tr }[c(\xi')b_0^{1}(x_0)]$ has no contribution for computing {\rm case~(c)}.

By (3.63) and (3.69), we have
\begin{eqnarray}
&&{\rm tr}\bigg\{\pi^+_{\xi_n}\sigma_{-1}(\widehat{D}^{-1}_{N})\times
\partial_{\xi_n}\bigg[\frac{c(\xi)}{|\xi|^6}
c(dx_n)\big[\partial_{x_n}[c(\xi')](x_0)|\xi|^2-c(\xi)h'(0)\big]+\frac{c(\xi)b_0^{2}(x_0)c(\xi)}{|\xi|^4}\bigg]
\bigg\}(x_0)|_{|\xi'|=1}\nonumber\\
&=&\frac{12h'(0)(i\xi^2_n+\xi_n-2i)}{(\xi-i)^3(\xi+i)^3}
+\frac{48h'(0)i\xi_n}{(\xi-i)^3(\xi+i)^4},
\end{eqnarray}
then
\begin{eqnarray}
-i\Omega_3\int_{\Gamma_+}\bigg[\frac{12h'(0)(i\xi_n^2+\xi_n-2i)}
{(\xi_n-i)^3(\xi_n+i)^3}+\frac{48h'(0)i\xi_n}{(\xi_n-i)^3(\xi_n+i)^4}\bigg]d\xi_ndx'=
-\frac{9}{2}\pi h'(0)\Omega_3dx'.
\end{eqnarray}
By (3.63) and (3.68), we have
\begin{eqnarray}
&&{\rm tr}\bigg[\pi^+_{\xi_n}\sigma_{-1}(\widehat{D}^{-1}_{N})\times
\partial_{\xi_n}\frac{c(\xi)(\bar{c}(\theta)-c(\theta'))c(\xi)}
{|\xi|^4}\bigg](x_0)|_{|\xi'|=1}\nonumber\\
&=&\frac{-1}{(\xi-i)(\xi+i)^3}{\rm tr}[c(\xi')(\bar{c}(\theta)-c(\theta'))(x_0)]
+\frac{i}{(\xi-i)(\xi+i)^3}{\rm tr}[c(dx_n)(\bar{c}(\theta)-c(\theta'))(x_0)].
\end{eqnarray}

By $\int_{|\xi'|=1}\{\xi_{1}\cdot\cdot\cdot\xi_{2d+1}\}\sigma(\xi')=0$ and (3.49), we have
\begin{eqnarray}
&&-i\int_{|\xi'|=1}\int^{+\infty}_{-\infty}{\rm tr}[\pi^+_{\xi_n}\sigma_{-1}(\widehat{D}^{-1}_{N})\times
\partial_{\xi_n}\frac{c(\xi)(\bar{c}(\theta)-c(\theta'))c(\xi)}
{|\xi|^4}](x_0)d\xi_n\sigma(\xi')dx'\nonumber\\
&=&-i\int_{|\xi'|=1}\int^{+\infty}_{-\infty}\frac{i}{(\xi-i)(\xi+i)^3}{\rm tr}[c(dx_n)(\bar{c}(\theta)-c(\theta'))](x_0)d\xi_n\sigma(\xi')dx'\nonumber\\
&=&-\frac{\pi}{4}{\rm tr}[c(dx_n)(\bar{c}(\theta)-c(\theta'))]\Omega_3dx'\nonumber\\
&=&-4\pi g(\theta',dx_n)\Omega_3dx'.
\end{eqnarray}
So we have
\begin{eqnarray}
{\rm case~ (c)}=-\frac{9}{2}\pi h'(0)\Omega_3dx'-4\pi g(\theta',dx_n)\Omega_3dx'.
\end{eqnarray}
Since $\Phi$ is the sum of the cases (a), (b) and (c), so
$$\Phi=-8\pi g(\theta',dx_n)\Omega_3dx'-2\pi i \Omega_3 f^{-1}\partial_{x_{n}}(f)dx'+2\pi i \Omega_3 f\partial_{x_{n}}(f^{-1})dx'.$$

\begin{thm}
Let $M$ be $4$-dimensional oriented
compact manifolds with the boundary $\partial M$ and the metric
$g^M$ as above, $\widehat{D}_{N}$ and $\widehat{D}^*_{N}$ be modified Novikov operators on $\widehat{M}$, then
\begin{eqnarray}
&&\widetilde{Wres}[\pi^{+}(f\widehat{D}^{-1}_{N}) \circ\pi^{+}\big(f^{-1}(\widehat{D}^{*}_{N})^{-1}\big)]\nonumber\\
&=&32\pi^2\int_{M}
\bigg\{\bigg[16g(\widetilde{e_{j}},\nabla^{TM}_{\widetilde{e_{j}}}\theta')-\frac{4}{3}s
-16|\theta|^2+32|\theta'|^2\bigg]+16\bigg[f^{-1}\Delta(f)+\langle grad_{M}f,grad_{M}f^{-1}\rangle\bigg]\nonumber\\
&&-16\Big(\bar{c}(\theta)-c(\theta')\Big)c(df)f^{-1}+8\bigg[-f^{-1}\Delta(f)(x_0)-\langle grad_{M}f,grad_{M}f^{-1}\rangle\bigg]+4\sum_{i}\Big[c(e_{i})c(\theta')\nonumber\\
&&-c(\theta')
c(e_{i})\Big]c(e_{i})c(df)f^{-1}+4\sum_{i}c(e_{i})c(df)f^{-1}
\Big[c(e_{i})c(\theta')-c(\theta')c(e_{i})\Big]  \bigg\}d{\rm Vol_{M}}\nonumber\\
&&+\int_{\partial M}\bigg\{-8\pi g(\theta',dx_n)\Omega_3dx'-2\pi i \Omega_3 f^{-1}\partial_{x_{n}}(f)dx'+2\pi i \Omega_3 f\partial_{x_{n}}(f^{-1})dx'\bigg\},
\end{eqnarray}
where $s$ is the scalar curvature.
\end{thm}

\section{ A Kastler-Kalau-Walze type theorem for six-dimensional
manifolds with boundary}

Let $M$ be $6$-dimensional compact manifolds with the boundary $\partial M$.
In the following, we will compute the more general case $\widetilde{Wres}[\pi^{+}(f\widehat{D}_{N}^{-1}) \circ\pi^{+}\big(f^{-1}(\widehat{D}_{N}^{*})^{-1}\cdot f\widehat{D}_{N}^{-1}\cdot f^{-1}(\widehat{D}^{*}_{N})^{-1}\big)]$ for nonzero
smooth functions $f,~f^{-1}$.
An application of (3.5) and (3.6) in \cite{Wa5} shows that
\begin{eqnarray}
&&\widetilde{Wres}[\pi^{+}(f\widehat{D}_{N}^{-1}) \circ\pi^{+}\big(f^{-1}(\widehat{D}_{N}^{*})^{-1}\cdot f\widehat{D}_{N}^{-1}\cdot f^{-1}(\widehat{D}^{*}_{N})^{-1}\big)]\nonumber\\
&=&\int_{M}\int_{|\xi|=1}{{\rm tr}}_{\wedge^*T^*M}\big[\sigma_{-n}\big( (\widehat{D}_{N}^{*}f\cdot\widehat{D}_{N}f^{-1})^{-2}\big)\big]\sigma(\xi)dx+\int_{\partial M}\Psi,
\end{eqnarray}
where
 \begin{eqnarray}
\Psi&=&\int_{|\xi'|=1}\int_{-\infty}^{+\infty}\sum_{j,k=0}^{\infty}\sum \frac{(-i)^{|\alpha|+j+k+\ell}}{\alpha!(j+k+1)!}
{{\rm tr}}_{\wedge^*T^*M}\Big[\partial_{x_{n}}^{j}\partial_{\xi'}^{\alpha}\partial_{\xi_{n}}^{k}\sigma_{r}^{+}
(f\widehat{D}_{N}^{-1})(x',0,\xi',\xi_{n})\nonumber\\
&&\times\partial_{x_{n}}^{\alpha}\partial_{\xi_{n}}^{j+1}\partial_{x_{n}}^{k}\sigma_{l}
\Big(f^{-1}(\widehat{D}_{N}^{*})^{-1}\cdot f\widehat{D}_{N}^{-1}\cdot f^{-1}(\widehat{D}^{*}_{N})^{-1}\Big)(x',0,\xi',\xi_{n})\Big]
d\xi_{n}\sigma(\xi')dx' ,
\end{eqnarray}
and the sum is taken over $r-k+|\alpha|+\ell-j-1=-n=-6,r\leq-1,\ell\leq-3$.

Note that
\begin{eqnarray}
&&f^{-1}(\widehat{D}_{N}^{*})^{-1}\cdot f\widehat{D}_{N}^{-1}\cdot f^{-1}(\widehat{D}_{N}^{*})^{-1}\nonumber\\
&=&(\widehat{D}_{N}^{*}f\cdot
\widehat{D}_{N}f^{-1}\cdot\widehat{D}_{N}^{*}f)^{-1}\nonumber\\
&=&\Big(\widehat{D}_{N}^{*}f\cdot
\widehat{D}_{N}\widehat{D}_{N}^{*}f^{-1}\cdot f-\widehat{D}_{N}^{*}f\cdot
\widehat{D}_{N}\cdot[\widehat{D}_{N}^{*},f^{-1}]\cdot f\Big)^{-1}\nonumber\\
&=&\Big(\widehat{D}_{N}^{*}f\cdot
\widehat{D}_{N}\widehat{D}_{N}^{*}-\widehat{D}_{N}^{*}f\cdot
\widehat{D}_{N}c(df^{-1})f\Big)^{-1}\nonumber\\
&=&\Big(f\cdot\widehat{D}_{N}^{*}\widehat{D}_{N}\widehat{D}_{N}^{*}
+[\widehat{D}_{N}^{*},f]\widehat{D}_{N}\widehat{D}_{N}^{*}
-\widehat{D}_{N}^{*}f\cdot\widehat{D}_{N}c(df^{-1})f\Big)^{-1}\nonumber\\
&=&\Big(f\cdot\widehat{D}_{N}^{*}\widehat{D}_{N}\widehat{D}_{N}^{*}
+c(df)\widehat{D}_{N}\widehat{D}_{N}^{*}
-\widehat{D}_{N}^{*}f\cdot\widehat{D}_{N}c(df^{-1})f\Big)^{-1}\nonumber\\
&=&\Big(f\cdot\widehat{D}_{N}^{*}\widehat{D}_{N}\widehat{D}_{N}^{*}
+c(df)\widehat{D}_{N}\widehat{D}_{N}^{*}
-\widehat{D}_{N}^{*}\widehat{D}_{N}f\cdot c(df^{-1})\cdot f
+\widehat{D}_{N}^{*}\cdot c(df) c(df^{-1})f\Big)^{-1}.
\end{eqnarray}
In order to get the symbol of operators $\widehat{D}_{N}^{*}f\cdot
\widehat{D}_{N}f^{-1}\cdot\widehat{D}_{N}^{*}f$. We first give the specification of
$\widehat{D}_{N}^{*}\widehat{D}_{N}\widehat{D}_{N}^{*}$, $\widehat{D}_{N}^{*}\widehat{D}_{N}$ and $\widehat{D}_{N}\widehat{D}_{N}^{*}$.
By (2.10) and (2.11), we have
\begin{eqnarray}
&&\widehat{D}_{N}\widehat{D}_{N}^{*}\nonumber\\
                &=&-\sum_{i,j}g^{i,j}\Big[\partial_{i}\partial_{j}
+2\sigma_{i}\partial_{j}+2a_{i}\partial_{j}-\Gamma_{i,j}^{k}\partial_{k}
+(\partial_{i}\sigma_{j})
+(\partial_{i}a_{j})
+\sigma_{i}\sigma_{j}+\sigma_{i}a_{j}+a_{i}\sigma_{j}+a_{i}a_{j} -\Gamma_{i,j}^{k}\sigma_{k}\nonumber\\
&&-\Gamma_{i,j}^{k}a_{k}\Big]
-\sum_{i,j}g^{i,j}\Big[c(\partial_{i})c(\theta')-c(\theta')c(\partial_{i})\Big]\partial_{j}
+\sum_{i,j}g^{i,j}\Big[c(\theta')c(\partial_{i})\sigma_{i}+c(\theta')c(\partial_{i})a_{i}
-c(\partial_{i})\partial_{i}(c(\theta'))\nonumber\\
&&-c(\partial_{i})\sigma_{i}c(\theta')-c(\partial_{i})a_{i}c(\theta')\Big]
-\frac{1}{8}\sum_{ijkl}R_{ijkl}\bar{c}(\widetilde{e_i})\bar{c}(\widetilde{e_j})
c(\widetilde{e_k})c(\widetilde{e_l})+\frac{1}{4}s
+\sum_{i}c(\widetilde{e_i})\bar{c}(\nabla_{\widetilde{e_i}}^{TM}\theta)+|\theta|^{2}\nonumber\\
&&+|\theta'|^2-\bar{c}(\theta)c(\theta')+c(\theta')\bar{c}(\theta)
\end{eqnarray}
and
\begin{eqnarray}
&&\widehat{D}^*_{N}\widehat{D}_{N}\nonumber\\
&=&-\sum_{i,j}g^{i,j}\Big[\partial_{i}\partial_{j}
+2\sigma_{i}\partial_{j}+2a_{i}\partial_{j}-\Gamma_{i,j}^{k}\partial_{k}
+(\partial_{i}\sigma_{j})
+(\partial_{i}a_{j})
+\sigma_{i}\sigma_{j}+\sigma_{i}a_{j}+a_{i}\sigma_{j}+a_{i}a_{j} -\Gamma_{i,j}^{k}\sigma_{k}\nonumber\\
&&-\Gamma_{i,j}^{k}a_{k}\Big]
+\sum_{i,j}g^{i,j}\Big[c(\partial_{i})c(\theta')-c(\theta')c(\partial_{i})\Big]\partial_{j}
-\sum_{i,j}g^{i,j}\Big[c(\theta')c(\partial_{i})\sigma_{i}+c(\theta')c(\partial_{i})a_{i}
-c(\partial_{i})\partial_{i}(c(\theta'))\nonumber\\
&&-c(\partial_{i})\sigma_{i}c(\theta')-c(\partial_{i})a_{i}c(\theta')\Big]
-\frac{1}{8}\sum_{ijkl}R_{ijkl}\bar{c}(\widetilde{e_i})\bar{c}(\widetilde{e_j})
c(\widetilde{e_k})c(\widetilde{e_l})+\frac{1}{4}s
+\sum_{i}c(\widetilde{e_i})\bar{c}(\nabla_{\widetilde{e_i}}^{TM}\theta)+|\theta|^{2}\nonumber\\
&&+|\theta'|^2+\bar{c}(\theta)c(\theta')-c(\theta')\bar{c}(\theta).
\end{eqnarray}
Combining (2.11) and (4.4), we obtain
\begin{eqnarray}
&&\widehat{D}^*_{N}\widehat{D}_{N}\widehat{D}^*_{N}\nonumber\\
&=&\sum^{n}_{i=1}c(\widetilde{e_{i}})\langle \widetilde{e_{i}},dx_{l}\rangle(-g^{ij}\partial_{l}\partial_{i}\partial_{j})
+\sum^{n}_{i=1}c(\widetilde{e_{i}})\langle \widetilde{e_{i}},dx_{l}\rangle \bigg\{-(\partial_{l}g^{ij})\partial_{i}\partial_{j}-g^{ij}\bigg(4(\sigma_{i}
+a_{i})\partial_{j}-2\Gamma^{k}_{ij}\partial_{k}\bigg)\partial_{l}\bigg\} \nonumber\\
&&+\sum^{n}_{i=1}c(\widetilde{e_{i}})\langle \widetilde{e_{i}},dx_{l}\rangle \bigg\{-2(\partial_{l}g^{ij})(\sigma_{i}+a_i)\partial_{j}+g^{ij}
(\partial_{l}\Gamma^{k}_{ij})\partial_{k}-2g^{ij}[(\partial_{l}\sigma_{i})
+(\partial_{l}a_i)]\partial_{j}
         +(\partial_{l}g^{ij})\Gamma^{k}_{ij}\partial_{k}\nonumber\\
         &&+\sum_{j,k}\Big[\partial_{l}\Big(c(\theta')c(\widetilde{e_{j}})
         -c(\widetilde{e_{j}}) c(\theta')\Big)\Big]\langle \widetilde{e_{j}},dx^{k}\rangle\partial_{k}
         +\sum_{j,k}\Big(c(\theta')c(\widetilde{e_{j}})-c(\widetilde{e_{j}})c(\theta')\Big)\Big[\partial_{l}\langle \widetilde{e_{j}},dx^{k}\rangle\Big]\partial_{k} \bigg\}\nonumber\\
         &&+\sum^{n}_{i=1}c(\widetilde{e_{i}})\langle \widetilde{e_{i}},dx_{l}\rangle\partial_{l}\bigg\{-g^{ij}\Big[
         (\partial_{i}\sigma_{j})+(\partial_{i}a_{j})+\sigma_{i}\sigma_{j}+\sigma_{i}a_{j}
         +a_{i}\sigma_{j}+a_{i}a_{j}-\Gamma_{i,j}^{k}\sigma_{k}-\Gamma_{i,j}^{k}a_{k}\nonumber\\
         &&+\sum_{i,j}g^{i,j}\Big[c(\theta')c(\partial_{i})\sigma_{i}
         +c(\theta')c(\partial_{i})a_{i}-c(\partial_{i})\partial_{i}(c(\theta'))
         -c(\partial_{i})\sigma_{i}c(\theta')-c(\partial_{i})a_{i}c(\theta')\Big]+\frac{1}{4}s\nonumber\\
         &&-\frac{1}{8}\sum_{ijkl}R_{ijkl}\bar{c}(\widetilde{e_i})\bar{c}(\widetilde{e_j})
         c(\widetilde{e_k})c(\widetilde{e_l})
         +\sum_{i}c(\widetilde{e_i})\bar{c}(\nabla_{\widetilde{e_i}}^{TM}\theta)
         +|\theta|^{2}+|\theta'|^2-\bar{c}(\theta)c(\theta')+c(\theta')\bar{c}(\theta)\bigg\}\nonumber\\
         &&+\Big[(\sigma_{i}+a_{i})+(\bar{c}(\theta)-c(\theta'))\Big](-g^{ij}\partial_{i}\partial_{j})
         +\sum^{n}_{i=1}c(\widetilde{e_{i}})\langle \widetilde{e_{i}},dx_{l}\rangle \bigg\{2\sum_{j,k}\Big[c(\theta')c(\widetilde{e_{j}})-c(\widetilde{e_{j}})c(\theta')\Big]\nonumber\\
         &&\times\langle \widetilde{e_{i}},dx_{k}\rangle\bigg\}\partial_{l}\partial_{k}
         +\Big[(\sigma_{i}+a_{i})+(\bar{c}(\theta)-c(\theta'))\Big]
         \bigg\{-\sum_{i,j}g^{i,j}\Big[2\sigma_{i}\partial_{j}+2a_{i}\partial_{j}
         -\Gamma_{i,j}^{k}\partial_{k}+(\partial_{i}\sigma_{j})\nonumber\\
         &&+(\partial_{i}a_{j})+\sigma_{i}\sigma_{j}+\sigma_{i}a_{j}+a_{i}\sigma_{j}+a_{i}a_{j} -\Gamma_{i,j}^{k}\sigma_{k}-\Gamma_{i,j}^{k}a_{k}\Big]-\sum_{i,j}g^{i,j}\Big[c(\partial_{i})c(\theta')
         -c(\theta')c(\partial_{i})\Big]\partial_{j}\nonumber\\
         &&+\sum_{i,j}g^{i,j}\Big[c(\theta')c(\partial_{i})\sigma_{i}+c(\theta')c(\partial_{i})a_{i}
         -c(\partial_{i})\partial_{i}(c(\theta'))-c(\partial_{i})\sigma_{i}c(\theta')
         -c(\partial_{i})a_{i}c(\theta')\Big]+\frac{1}{4}s+|\theta'|^2\nonumber\\
         &&-\frac{1}{8}\sum_{ijkl}R_{ijkl}\bar{c}(\widetilde{e_i})\bar{c}(\widetilde{e_j})
         c(\widetilde{e_k})c(\widetilde{e_l})
         +\sum_{i}c(\widetilde{e_i})\bar{c}(\nabla_{\widetilde{e_i}}^{TM}\theta)
         +|\theta|^{2}-\bar{c}(\theta)c(\theta')+c(\theta')\bar{c}(\theta)\bigg\}.
\end{eqnarray}

Thus, using (4.3)-(4.6), we get the specification of
$\widehat{D}_{N}^{*}f\cdot
\widehat{D}_{N}f^{-1}\cdot\widehat{D}_{N}^{*}f$.
\begin{eqnarray*}
&&\widehat{D}_{N}^{*}f\cdot
\widehat{D}_{N}f^{-1}\cdot\widehat{D}_{N}^{*}f\nonumber\\
&=&f\cdot\widehat{D}_{N}^{*}\widehat{D}_{N}\widehat{D}_{N}^{*}
+c(df)\widehat{D}_{N}\widehat{D}_{N}^{*}
-\widehat{D}_{N}^{*}\widehat{D}_{N}f\cdot c(df^{-1})\cdot f
+\widehat{D}_{N}^{*}\cdot c(df) c(df^{-1})f\nonumber\\
&=&f\cdot\Bigg\{
\sum^{n}_{i=1}c(\widetilde{e_{i}})\langle \widetilde{e_{i}},dx_{l}\rangle(-g^{ij}\partial_{l}\partial_{i}\partial_{j})
+\sum^{n}_{i=1}c(\widetilde{e_{i}})\langle \widetilde{e_{i}},dx_{l}\rangle \bigg\{-(\partial_{l}g^{ij})\partial_{i}\partial_{j}-g^{ij}\bigg(4(\sigma_{i}
+a_{i})\partial_{j}-2\Gamma^{k}_{ij}\partial_{k}\bigg)\partial_{l}\bigg\} \nonumber\\
&&+\sum^{n}_{i=1}c(\widetilde{e_{i}})\langle \widetilde{e_{i}},dx_{l}\rangle \bigg\{-2(\partial_{l}g^{ij})(\sigma_{i}+a_i)\partial_{j}+g^{ij}
(\partial_{l}\Gamma^{k}_{ij})\partial_{k}-2g^{ij}[(\partial_{l}\sigma_{i})
+(\partial_{l}a_i)]\partial_{j}
         +(\partial_{l}g^{ij})\Gamma^{k}_{ij}\partial_{k}\nonumber\\
         &&+\sum_{j,k}\Big[\partial_{l}\Big(c(\theta')c(\widetilde{e_{j}})
         -c(\widetilde{e_{j}}) c(\theta')\Big)\Big]\langle \widetilde{e_{j}},dx^{k}\rangle\partial_{k}
         +\sum_{j,k}\Big(c(\theta')c(\widetilde{e_{j}})-c(\widetilde{e_{j}})c(\theta')\Big)\Big[\partial_{l}\langle \widetilde{e_{j}},dx^{k}\rangle\Big]\partial_{k} \bigg\}\nonumber\\
         &&+\sum^{n}_{i=1}c(\widetilde{e_{i}})\langle \widetilde{e_{i}},dx_{l}\rangle\partial_{l}\bigg\{-g^{ij}\Big[
         (\partial_{i}\sigma_{j})+(\partial_{i}a_{j})+\sigma_{i}\sigma_{j}+\sigma_{i}a_{j}
         +a_{i}\sigma_{j}+a_{i}a_{j}-\Gamma_{i,j}^{k}\sigma_{k}-\Gamma_{i,j}^{k}a_{k}\nonumber\\
    \end{eqnarray*}
\begin{eqnarray}
         &&+\sum_{i,j}g^{i,j}\Big[c(\theta')c(\partial_{i})\sigma_{i}
         +c(\theta')c(\partial_{i})a_{i}-c(\partial_{i})\partial_{i}(c(\theta'))
         -c(\partial_{i})\sigma_{i}c(\theta')-c(\partial_{i})a_{i}c(\theta')\Big]+\frac{1}{4}s\nonumber\\
         &&-\frac{1}{8}\sum_{ijkl}R_{ijkl}\bar{c}(\widetilde{e_i})\bar{c}(\widetilde{e_j})
         c(\widetilde{e_k})c(\widetilde{e_l})
         +\sum_{i}c(\widetilde{e_i})\bar{c}(\nabla_{\widetilde{e_i}}^{TM}\theta)
         +|\theta|^{2}+|\theta'|^2-\bar{c}(\theta)c(\theta')+c(\theta')\bar{c}(\theta)\bigg\}\nonumber\\
         &&+\Big[(\sigma_{i}+a_{i})+(\bar{c}(\theta)-c(\theta'))\Big](-g^{ij}\partial_{i}\partial_{j})
         +\sum^{n}_{i=1}c(\widetilde{e_{i}})\langle \widetilde{e_{i}},dx_{l}\rangle \bigg\{2\sum_{j,k}\Big[c(\theta')c(\widetilde{e_{j}})-c(\widetilde{e_{j}})c(\theta')\Big]\nonumber\\
         &&\times\langle \widetilde{e_{i}},dx_{k}\rangle\bigg\}\partial_{l}\partial_{k}
         +\Big[(\sigma_{i}+a_{i})+(\bar{c}(\theta)-c(\theta'))\Big]
         \bigg\{-\sum_{i,j}g^{i,j}\Big[2\sigma_{i}\partial_{j}+2a_{i}\partial_{j}
         -\Gamma_{i,j}^{k}\partial_{k}+(\partial_{i}\sigma_{j})\nonumber\\
         &&+(\partial_{i}a_{j})+\sigma_{i}\sigma_{j}+\sigma_{i}a_{j}+a_{i}\sigma_{j}+a_{i}a_{j} -\Gamma_{i,j}^{k}\sigma_{k}-\Gamma_{i,j}^{k}a_{k}\Big]-\sum_{i,j}g^{i,j}\Big[c(\partial_{i})c(\theta')
         -c(\theta')c(\partial_{i})\Big]\partial_{j}\nonumber\\
         &&+\sum_{i,j}g^{i,j}\Big[c(\theta')c(\partial_{i})\sigma_{i}+c(\theta')c(\partial_{i})a_{i}
         -c(\partial_{i})\partial_{i}(c(\theta'))-c(\partial_{i})\sigma_{i}c(\theta')
         -c(\partial_{i})a_{i}c(\theta')\Big]+\frac{1}{4}s+|\theta'|^2\nonumber\\
         &&-\frac{1}{8}\sum_{ijkl}R_{ijkl}\bar{c}(\widetilde{e_i})\bar{c}(\widetilde{e_j})
         c(\widetilde{e_k})c(\widetilde{e_l})
         +\sum_{i}c(\widetilde{e_i})\bar{c}(\nabla_{\widetilde{e_i}}^{TM}\theta)
         +|\theta|^{2}-\bar{c}(\theta)c(\theta')+c(\theta')\bar{c}(\theta)\bigg\}
         \Bigg\}
         +c(df)\nonumber\\
         &&\times\Bigg\{-\sum_{i,j}g^{i,j}\Big[\partial_{i}\partial_{j}
+2\sigma_{i}\partial_{j}+2a_{i}\partial_{j}-\Gamma_{i,j}^{k}\partial_{k}
+(\partial_{i}\sigma_{j})
+(\partial_{i}a_{j})
+\sigma_{i}\sigma_{j}+\sigma_{i}a_{j}+a_{i}\sigma_{j}+a_{i}a_{j} -\Gamma_{i,j}^{k}\sigma_{k}\nonumber\\
&&-\Gamma_{i,j}^{k}a_{k}\Big]
-\sum_{i,j}g^{i,j}\Big[c(\partial_{i})c(\theta')-c(\theta')c(\partial_{i})\Big]\partial_{j}
+\sum_{i,j}g^{i,j}\Big[c(\theta')c(\partial_{i})\sigma_{i}+c(\theta')c(\partial_{i})a_{i}
-c(\partial_{i})\partial_{i}(c(\theta'))\nonumber\\
&&-c(\partial_{i})\sigma_{i}c(\theta')-c(\partial_{i})a_{i}c(\theta')\Big]
-\frac{1}{8}\sum_{ijkl}R_{ijkl}\bar{c}(\widetilde{e_i})\bar{c}(\widetilde{e_j})
c(\widetilde{e_k})c(\widetilde{e_l})+\frac{1}{4}s
+\sum_{i}c(\widetilde{e_i})\bar{c}(\nabla_{\widetilde{e_i}}^{TM}\theta)+|\theta|^{2}\nonumber\\
&&+|\theta'|^2-\bar{c}(\theta)c(\theta')+c(\theta')\bar{c}(\theta)
               \Bigg\}
               -\Bigg\{-\sum_{i,j}g^{i,j}\Big[\partial_{i}\partial_{j}
+2\sigma_{i}\partial_{j}+2a_{i}\partial_{j}-\Gamma_{i,j}^{k}\partial_{k}
+(\partial_{i}\sigma_{j})
+(\partial_{i}a_{j})\nonumber\\
&&+\sigma_{i}\sigma_{j}+\sigma_{i}a_{j}+a_{i}\sigma_{j}+a_{i}a_{j} -\Gamma_{i,j}^{k}\sigma_{k}-\Gamma_{i,j}^{k}a_{k}\Big]
+\sum_{i,j}g^{i,j}\Big[c(\partial_{i})c(\theta')-c(\theta')c(\partial_{i})\Big]\partial_{j}
\nonumber\\
&&-\sum_{i,j}g^{i,j}\Big[c(\theta')c(\partial_{i})\sigma_{i}
+c(\theta')c(\partial_{i})a_{i}
-c(\partial_{i})\partial_{i}(c(\theta'))-c(\partial_{i})\sigma_{i}c(\theta')-c(\partial_{i})a_{i}c(\theta')\Big]\nonumber\\
&&
-\frac{1}{8}\sum_{ijkl}R_{ijkl}\bar{c}(\widetilde{e_i})\bar{c}(\widetilde{e_j})
c(\widetilde{e_k})c(\widetilde{e_l})+\frac{1}{4}s
+\sum_{i}c(\widetilde{e_i})\bar{c}(\nabla_{\widetilde{e_i}}^{TM}\theta)+|\theta|^{2}+|\theta'|^2+\bar{c}(\theta)c(\theta')\nonumber\\
&&-c(\theta')\bar{c}(\theta)
               \Bigg\}f\cdot c(df^{-1})\cdot f+\Bigg\{\sum^{n}_{i,j=1}g^{ij}c(\partial_{i})\Big(\partial_{j}+a_{i}+\sigma_{i}\Big)+\bar{c}(\theta)-c(\theta')\Bigg\}\cdot c(df) c(df^{-1})f.
    \end{eqnarray}

Let $\partial^{j}=g^{ij}\partial_{i}$, 
by the above formulas, then we obtain:

\begin{lem}
Let $\widehat{D}_{N}$ and $\widehat{D}^*_{N}$ be modified Novikov operators on $\widehat{M}$, then
\begin{eqnarray}
\sigma_{3}(\widehat{D}_{N}^{*}f\cdot\widehat{D}_{N}f^{-1}\cdot\widehat{D}_{N}^{*}f)
&=& f\sigma_{3}(\widehat{D}_{N}^{*}\widehat{D}_{N}\widehat{D}_{N}^{*})
=\sqrt{-1}c(\xi)|\xi|^2f; \\
\sigma_{2}(\widehat{D}_{N}^{*}f\cdot\widehat{D}_{N}f^{-1}\cdot\widehat{D}_{N}^{*}f)
&=&f\sigma_{2}(\widehat{D}_{N}^{*}\widehat{D}_{N}\widehat{D}_{N}^{*})
+2c(df)|\xi|^2,
\end{eqnarray}
where
$\sigma_2(\widehat{D}^*_{N}\widehat{D}_{N}\widehat{D}^*_{N})=\sum_{i,j,l}c(dx_{l})\partial_{l}(g^{i,j})\xi_{i}\xi_{j} +c(\xi)(4\sigma^k+4a^k-2\Gamma^k)\xi_{k}-2\bigg[c(\xi)c(\theta')c(\xi)+|\xi|^2c(\theta')\bigg]$
$$+\frac{1}{4}|\xi|^2\sum_{s,t,l}\omega_{s,t}
(\widetilde{e_l})[c(\widetilde{e_l})\bar{c}(\widetilde{e_s})\bar{c}(\widetilde{e_t})
-c(\widetilde{e_l})c(\widetilde{e_s})c(\widetilde{e_t})]
+|\xi|^2(\bar{c}(\theta)-c(\theta')).$$

\end{lem}
In order to get the symbol of operators $\widehat{D}_{N}^{*}f\cdot
\widehat{D}_{N}f^{-1}\cdot\widehat{D}_{N}^{*}f$. We first give the following  formulas:
\begin{eqnarray}
D_x^{\alpha}=(-\sqrt{-1})^{|\alpha|}\partial_x^{\alpha};
~\sigma(\widehat{D}_{N}^{*}f\cdot\widehat{D}_{N}f^{-1}\cdot\widehat{D}_{N}^{*}f)&=&p_3+p_2+p_1+p_0;\nonumber\\
~\sigma\big((\widehat{D}_{N}^{*}f\cdot\widehat{D}_{N}f^{-1}
\cdot\widehat{D}_{N}^{*}f)^{-1}\big)&=&\sum^{\infty}_{j=3}q_{-j}.
\end{eqnarray}
By the composition formula of psudodifferential operators, we have
 \begin{eqnarray}
1
&=&\sigma\big[(\widehat{D}_{N}^{*}f\cdot\widehat{D}_{N}f^{-1}\cdot\widehat{D}_{N}^{*}f)\circ (\widehat{D}_{N}^{*}f\cdot\widehat{D}_{N}f^{-1}\cdot\widehat{D}_{N}^{*}f)^{-1}\big]\nonumber\\
&=&(p_3+p_2+p_1+p_0)(q_{-3}+q_{-4}+q_{-5}+\cdots) \nonumber\\
&&+\sum_j(\partial_{\xi_j}p_3+\partial_{\xi_j}p_2+\partial_{\xi_j}p_1+\partial_{\xi_j}p_0)
(D_{x_j}q_{-3}+D_{x_j}q_{-4}+D_{x_j}q_{-5}+\cdots) \nonumber\\
&=&p_3q_{-3}+(p_3q_{-4}+p_2q_{-3}+\sum_j\partial_{\xi_j}p_3D_{x_j}q_{-3})+\cdots.
\end{eqnarray}
Then
\begin{equation}
q_{-3}=p_3^{-1};~q_{-4}=-p_3^{-1}\Big[p_2p_3^{-1}+\sum_j\partial_{\xi_j}p_3D_{x_j}(p_3^{-1})\Big].
\end{equation}\\
By Lemma 2.1 in \cite{Wa3} and Lemma 4.1, we obtain

\begin{lem}Let $\widehat{D}_{N}$ and $\widehat{D}^*_{N}$ be modified Novikov operators on $\widehat{M}$, then
\begin{eqnarray}
\sigma_{-3}(\widehat{D}_{N}^{*}f\cdot\widehat{D}_{N}f^{-1}\cdot\widehat{D}_{N}^{*}f)^{-1}
&=&f^{-1}\sigma_{-3}(\widehat{D}_{N}^{*}\widehat{D}_{N}\widehat{D}_{N}^{*})^{-1}
=\frac{\sqrt{-1}c(\xi)}{f|\xi|^4}; \\
\sigma_{-4}(\widehat{D}_{N}^{*}f\cdot\widehat{D}_{N}f^{-1}\cdot\widehat{D}_{N}^{*}f)^{-1}
&=&f^{-1}\sigma_{-4}(\widehat{D}_{N}^{*}\widehat{D}_{N}\widehat{D}_{N}^{*})^{-1}
+\frac{2c(\xi)c(df)c(\xi)}{f^{2}|\xi|^6}\nonumber\\
&&+\frac{ic(\xi)\sum\limits_j\Big[c(dx_j)|\xi|^2+2\xi_{j}c(\xi)\Big]D_{x_j}(f^{-1})c(\xi)}{|\xi|^8},
\end{eqnarray}
where
\begin{eqnarray}
\sigma_{-3}\big((\widehat{D}^*_{N}\widehat{D}_{N}\widehat{D}^*_{N})^{-1}\big)&=&\frac{ic(\xi)}{|\xi|^{4}};\nonumber\\
\sigma_{-4}\big((\widehat{D}^*_{N}\widehat{D}_{N}\widehat{D}^*_{N})^{-1}\big)&=&
\frac{c(\xi)\sigma_{2}(\widehat{D}^*_{N}\widehat{D}_{N}\widehat{D}^*_{N})c(\xi)}{|\xi|^8}
+\frac{ic(\xi)}{|\xi|^8}\Big(|\xi|^4c(dx_n)\partial_{x_n}c(\xi')
-2h'(0)c(dx_n)c(\xi)\nonumber\\
&&+2\xi_{n}c(\xi)\partial_{x_n}c(\xi')+4\xi_{n}h'(0)\Big).
\end{eqnarray}
\end{lem}

Locally we can use Theorem 3.3 to compute the interior term of (4.1), then
 \begin{eqnarray}
&&\int_{M}\int_{|\xi|=1}{{\rm tr}}_{\wedge^*T^*M}[\sigma_{-n}\big( (\widehat{D}_{N}^{*}f\cdot\widehat{D}_{N}f^{-1})^{-2}\big)\big]\sigma(\xi)dx\nonumber\\
&=&128\pi^{3}\int_{M}
2^6\bigg\{\bigg(-\frac{1}{12}s-|\theta|^2+(n-2)|\theta'|^2+g(\widetilde{e_{j}},\nabla^{TM}_{\widetilde{e_{j}}}\theta')\bigg)+\bigg[f^{-1}\Delta(f)
+\langle grad_{M}f,grad_{M}f^{-1}\rangle\bigg]\nonumber\\
&&-\Big(\bar{c}(\theta)-c(\theta')\Big)c(df)f^{-1}
-\frac{1}{2}\bigg[f^{-1}\Delta(f)(x_0)+\langle grad_{M}f,grad_{M}f^{-1}\rangle\bigg]+\frac{1}{4}\sum_{i}\Big[c(e_{i})c(\theta')-c(\theta')
c(e_{i})\Big]\nonumber\\
&&\times c(e_{i})c(df)f^{-1}+\frac{1}{4}\sum_{i}c(e_{i})c(df)f^{-1}
\Big[c(e_{i})c(\theta')-c(\theta')c(e_{i})\Big]\bigg\}d{\rm Vol_{M}}.
\end{eqnarray}
So we only need to compute $\int_{\partial M}\Psi$.

From the formula (4.2) for the definition of $\Psi$, now we can compute $\Psi$.
Since the sum is taken over $r+\ell-k-j-|\alpha|-1=-6, \ r\leq-1, \ell\leq -3$, then we have the $\int_{\partial{M}}\Psi$
is the sum of the following five cases:
~\\
~\\
\noindent  {\bf case (a)~(I)}~$r=-1, l=-3, j=k=0, |\alpha|=1$.\\

By (4.2), we get
 \begin{eqnarray}
{\rm case~(a)~(I)}&=&-\int_{|\xi'|=1}\int^{+\infty}_{-\infty}\sum_{|\alpha|=1}{\rm tr}
\Big[\partial^{\alpha}_{\xi'}\pi^{+}_{\xi_{n}}\sigma_{-1}(f\widehat{D}_{N}^{-1})
      \times\partial^{\alpha}_{x'}\partial_{\xi_{n}}\sigma_{-3}
      \big(f^{-1}(\widehat{D}_{N}^{*})^{-1}\cdot f\widehat{D}_{N}^{-1}\cdot\nonumber\\
      &&f^{-1}(\widehat{D}_{N}^{*})^{-1}\big)\Big](x_0)d\xi_n\sigma(\xi')dx'\nonumber\\
&=&-\int_{|\xi'|=1}\int^{+\infty}_{-\infty}\sum_{|\alpha|=1}{\rm tr}
\Big[\partial^{\alpha}_{\xi'}\pi^{+}_{\xi_{n}}\big(f\sigma_{-1}(\widehat{D}_{N}^{-1})\big)
      \times\partial^{\alpha}_{x'}\partial_{\xi_{n}}\big(f^{-1}\sigma_{-3}
(\widehat{D}_{N}^{*}\widehat{D}_{N}\widehat{D}_{N}^{*})^{-1}\big)\Big](x_0)\nonumber\\
      &&\times d\xi_n\sigma(\xi')dx'\nonumber\\
  &=&-\int_{|\xi'|=1}\int^{+\infty}_{-\infty}\sum_{|\alpha|=1}{\rm tr}
\Big[\partial^{\alpha}_{\xi'}\pi^{+}_{\xi_{n}}\sigma_{-1}(\widehat{D}_{N}^{-1})
      \times\partial^{\alpha}_{x'}\partial_{\xi_{n}}\sigma_{-3}
(\widehat{D}_{N}^{*}\widehat{D}_{N}\widehat{D}_{N}^{*})^{-1}\Big](x_0)d\xi_n\sigma(\xi')dx'\nonumber\\
      &&-f\sum\limits_{j<n}\partial_{j}(f^{-1})\int_{|\xi'|=1}\int^{+\infty}_{-\infty}\sum_{|\alpha|=1}{\rm tr}
\Big[\partial^{\alpha}_{\xi'}\pi^{+}_{\xi_{n}}\sigma_{-1}(\widehat{D}_{N}^{-1})
      \times\partial_{\xi_{n}}\sigma_{-3}
(\widehat{D}_{N}^{*}\widehat{D}_{N}\widehat{D}_{N}^{*})^{-1}\Big](x_0)\nonumber\\
      &&\times d\xi_n\sigma(\xi')dx'.
\end{eqnarray}
By Lemma 2.2 in \cite{Wa3} and (4.15), for $i<n$, we have
 \begin{eqnarray}
 &&\partial_{x_{i}}\sigma_{-3}\Big((\widehat{D}_{N}^{*}\widehat{D}_{N}\widehat{D}_{N}^{*})^{-1}\Big)(x_0)=
      \partial_{x_{i}}\Big[\frac{\sqrt{-1}c(\xi)}{|\xi|^{4}}\Big](x_{0})\nonumber\\
      &=&\sqrt{-1}\partial_{x_{i}}\Big[c(\xi)\Big]|\xi|^{-4}(x_{0})
      -2\sqrt{-1}c(\xi)\partial_{x_{i}}\Big[|\xi|^{2}\Big]|\xi|^{-6}(x_{0})=0.
\end{eqnarray}
Thus we have
\begin{eqnarray}
-\int_{|\xi'|=1}\int^{+\infty}_{-\infty}\sum_{|\alpha|=1}{\rm tr}
\Big[\partial^{\alpha}_{\xi'}\pi^{+}_{\xi_{n}}\sigma_{-1}(\widehat{D}_{N}^{-1})
      \times\partial^{\alpha}_{x'}\partial_{\xi_{n}}\sigma_{-3}
(\widehat{D}_{N}^{*}\widehat{D}_{N}\widehat{D}_{N}^{*})^{-1}\Big](x_0)d\xi_n\sigma(\xi')dx'=0.
\end{eqnarray}
By (3.63) and direct calculations, for $i<n$, we obtain
\begin{eqnarray}
&&\partial^{\alpha}_{\xi'}\pi^{+}_{\xi_{n}}\sigma_{-1}(\widehat{D}_{N}^{-1})(x_0)|_{|\xi'|=1}
=\partial_{\xi_i}\pi^{+}_{\xi_{n}}\sigma_{-1}(\widehat{D}_{N}^{-1})(x_0)|_{|\xi'|=1}\nonumber\\
&=&\frac{c(dx_i)}{2(\xi_n-\sqrt{-1})}-\frac{\xi_i(\xi_n-2\sqrt{-1})c(\xi')+\xi_ic(dx_n)}{2(\xi_n-\sqrt{-1})^2},
\end{eqnarray}
and we get
\begin{eqnarray}
\partial_{\xi_{n}}\sigma_{-3}(\widehat{D}_{N}^{*}\widehat{D}_{N}\widehat{D}_{N}^{*})^{-1}
=\frac{\sqrt{-1}c(dx_n)}{|\xi|^4}-\frac{4\sqrt{-1}\big[\xi_nc(\xi')+\xi^2_nc(dx_n)\big]}{|\xi|^6}.
\end{eqnarray}
Then for $i<n$, we have
\begin{eqnarray}
&&{\rm tr}
\Big[\partial^{\alpha}_{\xi'}\pi^{+}_{\xi_{n}}\sigma_{-1}(\widehat{D}_{N}^{-1})
      \times\partial_{\xi_{n}}\sigma_{-3}
(\widehat{D}_{N}^{*}\widehat{D}_{N}\widehat{D}_{N}^{*})^{-1}\Big](x_0)\nonumber\\
&=&-\xi_i{\rm tr}
\Big[\frac{c(dx_n)^{2}}{2(\xi_n-\sqrt{-1})^2}\Big]-4\sqrt{-1}\xi_n\xi_i{\rm tr}
\Big[\frac{c(dx_i)^{2}}{2(\xi_n-\sqrt{-1})|\xi|^6}\Big]+4\sqrt{-1}\xi_n\xi_i(\xi_n-2\sqrt{-1})\nonumber\\
&&\times{\rm tr}
\Big[\frac{c(\xi')^{2}}{2(\xi_n-\sqrt{-1})^2|\xi|^6}\Big]+4\sqrt{-1}\xi^{2}_n\xi_i{\rm tr}
\Big[\frac{c(dx_n)^{2}}{2(\xi_n-\sqrt{-1})^2|\xi|^6}\Big].
\end{eqnarray}
We note that $i<n,~\int_{|\xi'|=1}\xi_i\sigma(\xi')=0$,
so
\begin{eqnarray}
&&-f\sum\limits_{j<n}\partial_{j}(f^{-1})\int_{|\xi'|=1}\int^{+\infty}_{-\infty}\sum_{|\alpha|=1}{\rm tr}
\Big[\partial^{\alpha}_{\xi'}\pi^{+}_{\xi_{n}}\sigma_{-1}(\widehat{D}_{N}^{-1})
      \times\partial_{\xi_{n}}\sigma_{-3}
(\widehat{D}_{N}^{*}\widehat{D}_{N}\widehat{D}_{N}^{*})^{-1}\Big](x_0) d\xi_n\sigma(\xi')dx'\nonumber\\
&=&0.
\end{eqnarray}
Then we have ${\bf case~(a)~(I)}=0$.
~\\

\noindent  {\bf case (a)~(II)}~$r=-1, l=-3, |\alpha|=k=0, j=1$.\\

By (4.2), we have
  \begin{eqnarray}
{\rm case~(a)~(II)}&=&-\frac{1}{2}\int_{|\xi'|=1}\int^{+\infty}_{-\infty} {\rm
tr} \Big[\partial_{x_{n}}\pi^{+}_{\xi_{n}}\sigma_{-1}(f\widehat{D}_{N}^{-1})
\times\partial^{2}_{\xi_{n}}\sigma_{-3}
\big(f^{-1}(\widehat{D}_{N}^{*})^{-1}\cdot f\widehat{D}_{N}^{-1}\cdot\nonumber\\
&&f^{-1}(\widehat{D}_{N}^{*})^{-1}\big)\Big](x_0)d\xi_n\sigma(\xi')dx'\nonumber\\
  &=&-\frac{1}{2}\int_{|\xi'|=1}\int^{+\infty}_{-\infty}{\rm tr}
\Big[\partial_{x_{n}}\pi^{+}_{\xi_{n}}\sigma_{-1}(\widehat{D}_{N}^{-1})
      \times\partial^2_{\xi_{n}}\sigma_{-3}
(\widehat{D}_{N}^{*}\widehat{D}_{N}\widehat{D}_{N}^{*})^{-1}\Big](x_0)d\xi_n\sigma(\xi')dx'\nonumber\\
      &&-\frac{1}{2}f^{-1}\partial_{x_{n}}(f)
      \int_{|\xi'|=1}\int^{+\infty}_{-\infty}{\rm tr}
\Big[\pi^{+}_{\xi_{n}}\sigma_{-1}(\widehat{D}_{N}^{-1})
      \times\partial^2_{\xi_{n}}\sigma_{-3}
(\widehat{D}_{N}^{*}\widehat{D}_{N}\widehat{D}_{N}^{*})^{-1}\Big](x_0)\nonumber\\
      &&\times d\xi_n\sigma(\xi')dx'.
\end{eqnarray}
By (2.3), (2.4) and (3.26), we have
\begin{equation}
\pi^{+}_{\xi_{n}}\partial_{x_{n}}\sigma_{-1}(\widehat{D}_{N}^{-1})(x_{0})|_{|\xi'|=1}
=\frac{\partial_{x_{n}}[c(\xi')](x_{0})}{2(\xi_{n}-\sqrt{-1})}
+\sqrt{-1}h'(0)\bigg[\frac{\sqrt{-1}c(\xi')}{4(\xi_{n}-\sqrt{-1})}
+\frac{c(\xi')+\sqrt{-1}c(dx_{n})}{4(\xi_{n}-\sqrt{-1})^{2}}\bigg].
\end{equation}
By (4.15) and direct calculations, we have\\
\begin{equation}
\partial_{\xi_{n}}\sigma_{-3}((\widehat{D}_{N}^{*}\widehat{D}_{N}\widehat{D}_{N}^{*})^{-1})
=\frac{-4\sqrt{-1}\xi_nc(\xi')+\sqrt{-1}(1-3\xi_n^{2})c(dx_n)}{(1+\xi_n^{2})^3}
\end{equation}
and
\begin{equation}
\partial^{2}_{\xi_{n}}\sigma_{-3}((\widehat{D}_{N}^{*}\widehat{D}_{N}\widehat{D}_{N}^{*})^{-1})
=\sqrt{-1}\bigg[\frac{(20\xi^{2}_{n}-4)c(\xi')+
12(\xi^{3}_{n}-\xi_{n})c(dx_{n})}{(1+\xi_{n}^{2})^{4}}\bigg].
\end{equation}

Since $n=6$, ${\rm tr}_{S(TM)\otimes F}[-{\rm id}]=64$. By the relation of the Clifford action and ${\rm tr}PQ={\rm tr}QP$,  then
\begin{eqnarray}
&&{\rm tr}[c(\xi')c(dx_{n})]=0; \ {\rm tr}[c(dx_{n})^{2}]=-64;\
{\rm tr}[c(\xi')^{2}](x_{0})|_{|\xi'|=1}=-64;\nonumber\\
&&{\rm tr}[\partial_{x_{n}}[c(\xi')]c(dx_{n})]=0; \
{\rm tr}[\partial_{x_{n}}c(\xi')c(\xi')](x_{0})|_{|\xi'|=1}=-32h'(0).
\end{eqnarray}
By (4.25)-(4.28), we get
\begin{eqnarray}
&&{\rm
tr} \Big[\partial_{x_{n}}\pi^{+}_{\xi_{n}}\sigma_{-1}(\widehat{D}_{N}^{-1})
      \times\partial^{2}_{\xi_{n}}\sigma_{-3}((\widehat{D}_{N}^{*}\widehat{D}_{N}
      \widehat{D}_{N}^{*})^{-1})\Big](x_0)\nonumber\\
&=&\frac{64h'(0)(-1-3\xi_{n}\sqrt{-1}+5\xi^{2}_{n}+3\sqrt{-1}\xi^{3}_{n})}{(\xi_{n}-\sqrt{-1})^{6}
(\xi_{n}+\sqrt{-1})^{4}}.
\end{eqnarray}
Then we obtain

\begin{eqnarray}
&&-\frac{1}{2}\int_{|\xi'|=1}\int^{+\infty}_{-\infty}{\rm tr}
\Big[\partial_{x_{n}}\pi^{+}_{\xi_{n}}\sigma_{-1}(\widehat{D}_{N}^{-1})
      \times\partial^2_{\xi_{n}}\sigma_{-3}
(\widehat{D}_{N}^{*}\widehat{D}_{N}\widehat{D}_{N}^{*})^{-1}\Big](x_0)d\xi_n\sigma(\xi')dx'\nonumber\\
      &=&-\frac{15}{2}\pi h'(0)\Omega_{4}dx'.
\end{eqnarray}

On the other hand, by calculations, we have
\begin{equation}
\pi^{+}_{\xi_{n}}\sigma_{-1}(\widehat{D}_{N}^{-1})(x_{0})|_{|\xi'|=1}
=-\frac{c(\xi')+\sqrt{-1}c(dx_{n})}{2(\xi_{n}-\sqrt{-1})}.
\end{equation}

By (4.25), (4.27) and (4.31), we get
\begin{equation}
{\rm
tr} \Big[\pi^{+}_{\xi_{n}}\sigma_{-1}(\widehat{D}_{N}^{-1})
      \times\partial^{2}_{\xi_{n}}\sigma_{-3}((\widehat{D}_{N}^{*}
      \widehat{D}_{N}\widehat{D}_{N}^{*})^{-1})\Big](x_0)
=\frac{-128(5\xi^2_{n}\sqrt{-1}-\sqrt{-1}-3\xi^{3}_{n}+3\xi_{n})}
{(\xi_{n}-\sqrt{-1})^{5}(\xi_{n}+\sqrt{-1})^{4}}.
\end{equation}
Then we obtain
\begin{eqnarray}
&&-\frac{1}{2}f^{-1}\partial_{x_{n}}(f)
      \int_{|\xi'|=1}\int^{+\infty}_{-\infty}{\rm tr}
\Big[\pi^{+}_{\xi_{n}}\sigma_{-1}(\widehat{D}_{N}^{-1})
      \times\partial^2_{\xi_{n}}\sigma_{-3}
(\widehat{D}_{N}^{*}\widehat{D}_{N}\widehat{D}_{N}^{*})^{-1}\Big](x_0) d\xi_n\sigma(\xi')dx'\nonumber\\
      &=&\frac{5\sqrt{-1}+44}{4}\pi f^{-1}\partial_{x_{n}}(f)\cdot8{\rm \Omega_{4}}dx',
\end{eqnarray}
where ${\rm \Omega_{4}}$ is the canonical volume of $S_{4}.$

Combining (4.24), (4.30) and (4.33), we obtain
\begin{eqnarray}
{\bf case~(a)~(II)}&=&-\frac{15}{2}\pi h'(0)\Omega_{4}dx'+
(10\sqrt{-1}+88)\pi f^{-1}\partial_{x_{n}}(f)\cdot{\rm \Omega_{4}}dx'.
\end{eqnarray}

\noindent  {\bf case (a)~(III)}~$r=-1,l=-3,|\alpha|=j=0,k=1$.\\

By (4.2), we have
 \begin{eqnarray}
{\rm case~ (a)~(III)}&=&-\frac{1}{2}\int_{|\xi'|=1}\int^{+\infty}_{-\infty}{\rm tr} \Big[\partial_{\xi_{n}}\pi^{+}_{\xi_{n}}\sigma_{-1}(f\widehat{D}_{N}^{-1})
      \times\partial_{\xi_{n}}\partial_{x_{n}}
      \sigma_{-3}
\big(f^{-1}(\widehat{D}_{N}^{*})^{-1}\cdot f\widehat{D}_{N}^{-1}\cdot\nonumber\\
&&f^{-1}(\widehat{D}_{N}^{*})^{-1}\big)\Big](x_0)d\xi_n\sigma(\xi')dx'\nonumber\\
  &=&-\frac{1}{2}\int_{|\xi'|=1}\int^{+\infty}_{-\infty}{\rm tr}
\Big[\partial_{\xi_{n}}\pi^{+}_{\xi_{n}}\big(\sigma_{-1}(\widehat{D}_{N}^{-1})\big)
      \times\partial_{\xi_{n}}\partial_{x_{n}}\sigma_{-3}
(\widehat{D}_{N}^{*}\widehat{D}_{N}\widehat{D}_{N}^{*})^{-1}\Big](x_0)d\xi_n\sigma(\xi')dx'\nonumber\\
      &&-\frac{1}{2}f\partial_{x_{n}}(f^{-1})
      \int_{|\xi'|=1}\int^{+\infty}_{-\infty}{\rm tr}
\Big[\partial_{\xi_{n}}\pi^{+}_{\xi_{n}}\sigma_{-1}(\widehat{D}_{N}^{-1})
      \times\partial_{\xi_{n}}\sigma_{-3}
(\widehat{D}_{N}^{*}\widehat{D}_{N}\widehat{D}_{N}^{*})^{-1}\Big](x_0) d\xi_n\sigma(\xi')\nonumber\\
      &&\times dx'.
\end{eqnarray}
By (4.6) and direct calculations, we have\\
\begin{equation}
\partial_{\xi_{n}}\partial_{x_{n}}\sigma_{-3}((\widehat{D}_{N}^{*}\widehat{D}_{N}\widehat{D}_{N}^{*})^{-1})=-\frac{4 \sqrt{-1}\xi_{n}\partial_{x_{n}}c(\xi')(x_{0})}{(1+\xi_{n}^{2})^{3}}
      +\frac{12\sqrt{-1}h'(0)\xi_{n}c(\xi')}{(1+\xi_{n}^{2})^{4}}
      -\frac{\sqrt{-1}(2-10\xi^{2}_{n})h'(0)c(dx_{n})}{(1+\xi_{n}^{2})^{4}}.
\end{equation}

Combining (3.37) and (4.6), we have
\begin{eqnarray}
&&{\rm tr} \Big[\partial_{\xi_{n}}\pi^{+}_{\xi_{n}}\sigma_{-1}(\widehat{D}_{N}^{-1})
      \times\partial_{\xi_{n}}\partial_{x_{n}}\sigma_{-3}((\widehat{D}_{N}^{*}
      \widehat{D}_{N}\widehat{D}_{N}^{*})^{-1})\Big](x_{0})|_{|\xi'|=1}\nonumber\\
&=&\frac{64h'(0)(\sqrt{-1}-4\xi_{n}-\sqrt{-1}\xi^{2}_{n})}{(\xi_{n}-\sqrt{-1})^{5}(\xi+i)^{4}},
\end{eqnarray}
and
\begin{eqnarray}
&&{\rm tr} \Big[\partial_{\xi_{n}}\pi^{+}_{\xi_{n}}\sigma_{-1}(\widehat{D}_{N}^{-1})
      \times\partial_{\xi_{n}}\sigma_{-3}((\widehat{D}_{N}^{*}\widehat{D}_{N}
      \widehat{D}_{N}^{*})^{-1})\Big](x_{0})|_{|\xi'|=1}\nonumber\\
&=&-32\frac{4\sqrt{-1}\xi_n+1-3\xi^2_{n}}{(\xi_{n}-\sqrt{-1})^{5}(\xi_n+\sqrt{-1})^{3}}.
\end{eqnarray}
Then
\begin{eqnarray}
&&-\frac{1}{2}\int_{|\xi'|=1}\int^{+\infty}_{-\infty}{\rm tr}
\Big[\partial_{x_{n}}\pi^{+}_{\xi_{n}}\big(\sigma_{-1}(\widehat{D}_{N}^{-1})\big)
      \times\partial_{\xi_{n}}\partial_{x_{n}}\sigma_{-3}
(\widehat{D}_{N}^{*}\widehat{D}_{N}\widehat{D}_{N}^{*})^{-1}\Big](x_0)d\xi_n\sigma(\xi')dx'\nonumber\\
      &=&\frac{25}{2}\pi h'(0)\Omega_{4}dx',
\end{eqnarray}
and
\begin{eqnarray}
&&-\frac{1}{2}f\partial_{x_{n}}(f^{-1})
      \int_{|\xi'|=1}\int^{+\infty}_{-\infty}{\rm tr}
\Big[\partial_{\xi_{n}}\pi^{+}_{\xi_{n}}\sigma_{-1}(\widehat{D}_{N}^{-1})
      \times\partial_{\xi_{n}}\sigma_{-3}
(\widehat{D}_{N}^{*}\widehat{D}_{N}\widehat{D}_{N}^{*})^{-1}\Big](x_0) d\xi_n\sigma(\xi')dx'\nonumber\\
      &=&\frac{\pi \sqrt{-1}}{2}\cdot f\cdot\partial_{x_{n}}(f^{-1})\Omega_{4}dx',
\end{eqnarray}
where ${\rm \Omega_{4}}$ is the canonical volume of $S_{4}.$

Then
\begin{eqnarray}
{\bf case~(a)~(III)}&=&\Big[\frac{25}{2}\pi h'(0)+
\frac{\pi \sqrt{-1}}{2}\cdot f\cdot\partial_{x_{n}}(f^{-1})\Big]\Omega_{4}dx'.
\end{eqnarray}
\\
\noindent  {\bf case (b)}~$r=-1,l=-4,|\alpha|=j=k=0$.\\

By (4.2), we have
 \begin{eqnarray}
{\rm case~ (b)}&=&-i\int_{|\xi'|=1}\int^{+\infty}_{-\infty}{\rm tr} \Big[\pi^{+}_{\xi_{n}}\sigma_{-1}(f\widehat{D}_{N}^{-1})
      \times\partial_{\xi_{n}}\sigma_{-4}
      \big(f^{-1}(\widehat{D}_{N}^{*})^{-1}\cdot f\widehat{D}_{N}^{-1}\cdot f^{-1}(\widehat{D}_{N}^{*})^{-1}\big)\Big](x_0)\nonumber\\
&&\times d\xi_n\sigma(\xi')dx'\nonumber\\
&=&-i\int_{|\xi'|=1}\int^{+\infty}_{-\infty}{\rm tr} \Bigg[\pi^{+}_{\xi_{n}}\sigma_{-1}(f\widehat{D}_{N}^{-1})
      \times\partial_{\xi_{n}}\Bigg(f^{-1}\sigma_{-4}(\widehat{D}_{N}^{*}
      \widehat{D}_{N}\widehat{D}_{N}^{*})^{-1}+\frac{2c(\xi)c(df)c(\xi)}{f^{2}|\xi|^6}
\nonumber\\
      &&+\frac{ic(\xi)\sum\limits_j\Big[c(dx_j)|\xi|^2+2\xi_{j}c(\xi)\Big]D_{x_j}(f^{-1})c(\xi)}{|\xi|^8}
      \Bigg)\Bigg](x_0)d\xi_n\sigma(\xi')dx'\nonumber\\
     &=&-i\int_{|\xi'|=1}\int^{+\infty}_{-\infty}{\rm tr} \Big[\pi^{+}_{\xi_{n}}\sigma_{-1}(\widehat{D}_{N}^{-1})
      \times\partial_{\xi_{n}}\Big(\sigma_{-4}(\widehat{D}_{N}^{*}
      \widehat{D}_{N}\widehat{D}_{N}^{*})^{-1}\Big)\Big](x_0)d\xi_n\sigma(\xi')dx'\nonumber\\
&&-2if^{-1}\int_{|\xi'|=1}\int^{+\infty}_{-\infty}{\rm tr} \Bigg[\pi^{+}_{\xi_{n}}\sigma_{-1}(\widehat{D}_{N}^{-1})
      \times\partial_{\xi_{n}}\Bigg(\frac{c(\xi)c(df)c(\xi)}{|\xi|^6}
      \Bigg)\Bigg](x_0)d\xi_n\sigma(\xi')dx'\nonumber\\
      &&-fi\int_{|\xi'|=1}\int^{+\infty}_{-\infty}{\rm tr} \Big[\pi^{+}_{\xi_{n}}\sigma_{-1}(\widehat{D}_{N}^{-1})
      \times\partial_{\xi_{n}}\Big(\frac{ic(\xi)\sum\limits_j\Big[c(dx_j)|\xi|^2
      +2\xi_{j}c(\xi)\Big]D_{x_j}(f^{-1})c(\xi)}{|\xi|^8}\Big)\Big]\nonumber\\
      &&\times(x_0)d\xi_n\sigma(\xi')dx'.
\end{eqnarray}

In the normal coordinate, $g^{ij}(x_{0})=\delta^{j}_{i}$ and $\partial_{x_{j}}(g^{\alpha\beta})(x_{0})=0$, if $j<n$; $\partial_{x_{j}}(g^{\alpha\beta})(x_{0})=h'(0)\delta^{\alpha}_{\beta}$, if $j=n$.
So by Lemma A.2 in \cite{Wa3}, we have $\Gamma^{n}(x_{0})=\frac{5}{2}h'(0)$ and $\Gamma^{k}(x_{0})=0$ for $k<n$. By the definition of $\delta^{k}$ and Lemma 2.3 in \cite{Wa3}, we have $\delta^{n}(x_{0})=0$ and $\delta^{k}=\frac{1}{4}h'(0)c(\widetilde{e_{k}})c(\widetilde{e_{n}})$ for $k<n$. By Lemma 4.2, we obtain

\begin{eqnarray}
&&\sigma_{-4}((\widehat{D}_{N}^{*}\widehat{D}_{N}\widehat{D}_{N}^{*})^{-1})(x_{0})|_{|\xi'|=1}\nonumber\\
&=&
\frac{c(\xi)\sigma_{2}((\widehat{D}_{N}^{*}\widehat{D}_{N}\widehat{D}_{N}^{*})^{-1})
(x_{0})|_{|\xi'|=1}c(\xi)}{|\xi|^8}
-\frac{c(\xi)}{|\xi|^4}\sum_j\partial_{\xi_j}\big(c(\xi)|\xi|^2\big)
D_{x_j}\big(\frac{ic(\xi)}{|\xi|^4}\big)\nonumber\\
&=&\frac{1}{|\xi|^8}c(\xi)\Big(\frac{1}{2}h'(0)c(\xi)\sum_{k<n}\xi_k
c(\widetilde{e_k})c(\widetilde{e_n})-\frac{1}{2}h'(0)c(\xi)\sum_{k<n}\xi_k
\bar{c}(\widetilde{e_k})\bar{c}(\widetilde{e_n})-\frac{5}{2}h'(0)\xi_nc(\xi)\nonumber\\
&&-\frac{1}{4}h'(0)|\xi|^2c(dx_n)
-2[c(\xi)c(\theta')c(\xi)+|\xi|^2c(\theta')]+|\xi|^2(\bar{c}(\theta)-c(\theta'))\Big)c(\xi)\nonumber\\
&&
+\frac{ic(\xi)}{|\xi|^8}\Big(|\xi|^4c(dx_n)\partial_{x_n}c(\xi')
-2h'(0)c(dx_n)c(\xi)+2\xi_{n}c(\xi)\partial_{x_n}c(\xi')+4\xi_{n}h'(0)\Big).
\end{eqnarray}
By (3.37) and (4.43), we have
\begin{eqnarray}
&&{\rm tr} [\partial_{\xi_n}\pi^+_{\xi_n}\sigma_{-1}(\widehat{D}^{-1}_{N})\times
\sigma_{-4}(\widehat{D}_{N}^{*}\widehat{D}_{N}\widehat{D}_{N}^{*})^{-1}](x_0)|_{|\xi'|=1} \nonumber\\
&=&\frac{1}{2(\xi_{n}-i)^{2}(1+\xi_{n}^{2})^{4}}
\big[\frac{3}{4}i+2+(3+4i)\xi_{n}+(-6+2i)\xi_{n}^{2}+3\xi_{n}^{3}+\frac{9i}{4}\xi_{n}^{4}\big]h'(0){\rm tr}
[{\rm id}]\nonumber\\
&&+\frac{1}{2(\xi_{n}-i)^{2}(1+\xi_{n}^{2})^{4}}\big(-1-3i\xi_{n}-2\xi_{n}^{2}-4i\xi_{n}^{3}-\xi_{n}^{4}-i\xi_{n}^{5}\big){\rm tr[c(\xi')\partial_{x_n}c(\xi')]}\nonumber\\
&&-\frac{1}{2(\xi_{n}-i)^{2}(1+\xi_{n}^{2})^{4}}\big(\frac{1}{2}i+\frac{1}{2}\xi_{n}+\frac{1}{2}\xi_{n}^{2}+\frac{1}{2}\xi_{n}^{3}\big){\rm tr}
[c(\xi')\bar{c}(\xi')c(dx_n)\bar{c}(dx_n)]\nonumber\\
&&+\frac{-\xi_ni+3}{2(\xi_{n}-i)^{4}(i+\xi_{n})^{3}}{\rm tr}\big[c(\theta')c(dx_n)\big]-\frac{3\xi_n+i}{2(\xi_{n}-i)^{4}(i+\xi_{n})^{3}}{\rm tr}\big[c(\theta')c(\xi')\big].
\end{eqnarray}
By direct calculation and the relation of the Clifford action and ${\rm tr}{AB}={\rm tr }{BA}$, then we have equalities:
\begin{eqnarray}
&&{\rm tr }[c(\theta')(x_0)c(dx_n)]=-64g(\theta',dx_n);~~{\rm tr }[c(\theta')(x_0)c(\xi')]=-64g(\theta',\xi');\nonumber\\
&&{\rm tr}[c(\widetilde{e_i})
\bar{c}(\widetilde{e_i})c(\widetilde{e_n})
\bar{c}(\widetilde{e_n})]=0~~(i<n).
\end{eqnarray}
Then
\begin{eqnarray}
{\rm tr}
[c(\xi')\bar{c}(\xi')c(dx_n)\bar{c}(dx_n)]&=&
\sum_{i<n,j<n}{\rm tr}[\xi_{i}\xi_{j}c(\widetilde{e_i})\bar{c}
(\widetilde{e_j})c(dx_n)\bar{c}(dx_n)]=0.
\end{eqnarray}
So, we have
\begin{eqnarray}
&&-i\int_{|\xi'|=1}\int^{+\infty}_{-\infty}{\rm tr} \Big[\pi^{+}_{\xi_{n}}\sigma_{-1}(\widehat{D}_{N}^{-1})
      \times\partial_{\xi_{n}}\Big(\sigma_{-4}(\widehat{D}^{*}_{N}
      \widehat{D}_{N}\widehat{D}^{*}_{N})^{-1}\Big)\Big](x_0)d\xi_n\sigma(\xi')dx'\nonumber\\
      &=&
 ih'(0)\int_{|\xi'|=1}\int^{+\infty}_{-\infty}64\times\frac{\frac{3}{4}i+2+(3+4i)\xi_{n}+(-6+2i)\xi_{n}^{2}+3\xi_{n}^{3}+\frac{9i}{4}\xi_{n}^{4}}{2(\xi_n-i)^5(\xi_n+i)^4}d\xi_n\sigma(\xi')dx'\nonumber\\ &&+ih'(0)\int_{|\xi'|=1}\int^{+\infty}_{-\infty}32\times\frac{1+3i\xi_{n}
 +2\xi_{n}^{2}+4i\xi_{n}^{3}+\xi_{n}^{4}+i\xi_{n}^{5}}{2(\xi_{n}-i)^{2}
 (1+\xi_{n}^{2})^{4}}d\xi_n\sigma(\xi')dx'\nonumber\\
 &&+i\int_{|\xi'|=1}\int^{+\infty}_{-\infty}
 \frac{\xi_n-i-2\xi_n i +1}{2(\xi_{n}-i)^{4}(i+\xi_{n})^{3}}{\rm tr}\big[c(\theta')c(dx_n)\big]d\xi_n\sigma(\xi')dx'\nonumber\\
 &&-i\int_{|\xi'|=1}\int^{+\infty}_{-\infty}
 \frac{3\xi_n+i}{2(\xi_{n}-i)^{4}(i+\xi_{n})^{3}}{\rm tr}\big[c(\theta')c(\xi')\big]d\xi_n\sigma(\xi')dx'\nonumber\\
&=&(-\frac{19}{4}i-15)\pi h'(0)\Omega_4dx'+(-\frac{3}{8}i-\frac{75}{8})\pi h'(0)\Omega_4dx'+120i\pi g(dx_n,\theta')\Omega_4dx'\nonumber\\
&=&(-\frac{41}{8}i-\frac{195}{8})\pi h'(0)\Omega_4dx'+120i\pi g(dx_n,\theta')\Omega_4dx'.
\end{eqnarray}

Since
\begin{eqnarray*}
\partial_{\xi_{n}}\Big(\frac{c(\xi)c(df) c(\xi)}{|\xi|^{6}}\Big)&=&\frac{c(dx_{n})c(df) c(\xi')+c(\xi')c(df) c(dx_{n})+2\xi_{n}c(dx_{n})c(df) c(dx_{n})}{(1+\xi^{2}_{n})^{3}}\nonumber\\
&&-\frac{6\xi_{n}c(\xi)c(df) c(\xi)}{(1+\xi^{2}_{n})^{4}}
\end{eqnarray*}
and
\begin{eqnarray}
&&\partial_{\xi_{n}}\Big(\frac{ic(\xi)\sum\limits_j\Big[c(dx_j)|\xi|^2
+2\xi_{j}c(\xi)\Big]D_{x_j}(f^{-1})c(\xi)}{|\xi|^8}\Big)\nonumber\\
&=&i\Bigg\{c(dx_{n})\sum\limits_j\Big[c(dx_j)|\xi|^2
+2\xi_{j}c(\xi)\Big]D_{x_j}(f^{-1}) c(\xi')
+c(\xi')\sum\limits_j\Big[c(dx_j)|\xi|^2
+2\xi_{j}c(\xi)\Big]D_{x_j}(f^{-1}) c(dx_{n})\nonumber\\
&&+2\xi_{n}c(dx_{n})\sum\limits_j\Big[c(dx_j)|\xi|^2
+2\xi_{j}c(\xi)\Big]D_{x_j}(f^{-1}) c(dx_{n})\Bigg\}(1+\xi^{2}_{n})^{-4}-i\Bigg\{8\xi_{n}c(\xi)\sum\limits_j\Big[c(dx_j)|\xi|^2
\nonumber\\
&&+2\xi_{j}c(\xi)\Big]D_{x_j}(f^{-1}) c(\xi)\Bigg\}(1+\xi^{2}_{n})^{-5},
\end{eqnarray}
then we have
\begin{eqnarray*}
&&{\rm tr} \Big[\pi^{+}_{\xi_{n}}\sigma_{-1}(\widehat{D}_{N}^{-1})
      \times\partial_{\xi_{n}}\Big(\frac{c(\xi)c(df) c(\xi)}{|\xi|^{6}}\Big)\Big](x_0) \nonumber\\
&=&\frac{(4\xi_{n}i+2)i}{2(\xi_{n}+i)(1+\xi^{2}_{n})^{3}}{\rm tr}[c(\xi')c(df)]+\frac{4\xi_{n}i+2}{2(\xi_{n}+i)(1+\xi^{2}_{n})^{3}}{\rm tr}[c(dx_{n})c(df)]
\end{eqnarray*}
and
\begin{eqnarray}
&&{\rm tr} \Bigg[\pi^{+}_{\xi_{n}}\sigma_{-1}(\widehat{D}_{N}^{-1})
      \times\partial_{\xi_{n}}\Big(\frac{ic(\xi)\sum\limits_j\Big[c(dx_j)|\xi|^2
+2\xi_{j}c(\xi)\Big]D_{x_j}(f^{-1})c(\xi)}{|\xi|^8}\Big)\Bigg](x_0) \nonumber\\
&=&\frac{(3\xi_{n}-i)i}{(\xi_{n}+i)(1+\xi^{2}_{n})^{4}}{\rm tr}\Bigg[c(\xi')\sum\limits_j\Big[c(dx_j)|\xi|^2
+2\xi_{j}c(\xi)\Big]D_{x_j}(f^{-1})\Bigg]\nonumber\\
&&+\frac{3\xi_{n}-i}{(\xi_{n}+i)(1+\xi^{2}_{n})^{4}}{\rm tr}\Bigg[c(dx_{n})\sum\limits_j\Big[c(dx_j)|\xi|^2
+2\xi_{j}c(\xi)\Big]D_{x_j}(f^{-1})\Bigg].
\end{eqnarray}

By the relation of the Clifford action and ${\rm tr}QP={\rm tr}PQ$, then we have the following equalities
\begin{eqnarray*}
&&{\rm tr}\Big[c(dx_{n})c(df)\Big]=-g(dx_{n},df)
\end{eqnarray*}
and
\begin{eqnarray*}
&&{\rm tr}\Bigg[c(dx_{n})\sum\limits_j\Big[c(dx_j)|\xi|^2
      +2\xi_{j}c(\xi)\Big]D_{x_j}(f^{-1})\Bigg]\nonumber\\
      &=&{\rm tr}\big(-{\rm id}\big)|\xi|^2\bigg(-i\partial_{x_n}(f)f^{-1}\bigg)
      +2\sum\limits_{j}\xi_j\xi_n{\rm tr}\big(-{\rm id}\big)\bigg(-i\partial_{x_j}(f)f^{-1}\bigg)\nonumber\\
      &=&-64|\xi|^2\bigg(-i\partial_{x_n}(f)f^{-1}\bigg)+2\sum\limits_{j}\xi_j\xi_n{\rm tr}\big(-{\rm id}\big)\bigg(-i\partial_{x_j}(f)f^{-1}\bigg).
\end{eqnarray*}

We note that $i<n,~\int_{|\xi'|=1}\xi_i\sigma(\xi')=0$,
so ${\rm tr}\big[c(\xi')c(df)\big]$, ${\rm tr}\Bigg[c(\xi')\sum\limits_j\Big[c(dx_j)|\xi|^2
      +2\xi_{j}c(\xi)\Big]D_{x_j}(f^{-1})\Bigg]$ and $2i\sum\limits_{j}\xi_j\xi_n\partial_{x_j}(f)f^{-1}{\rm tr}[-{\rm id}]$ have no contribution for computing {\bf case (b)}.
Then we obtain
\begin{eqnarray}
&&-2if^{-1}\int_{|\xi'|=1}\int^{+\infty}_{-\infty}{\rm tr} \Bigg[\pi^{+}_{\xi_{n}}\sigma_{-1}(\widehat{D}_{N}^{-1})
      \times\partial_{\xi_{n}}\Big(\frac{c(\xi)c(df) c(\xi)}{|\xi|^{6}}\Big)\Bigg](x_0)d\xi_n\sigma(\xi')dx'\nonumber\\
      &=&\frac{3}{8f}\pi  g(dx_{n},df)\Omega_{4}dx'.
\end{eqnarray}
and
\begin{eqnarray}
&&-fi\int_{|\xi'|=1}\int^{+\infty}_{-\infty}{\rm tr} \Big[\pi^{+}_{\xi_{n}}\sigma_{-1}(\widehat{D}_{N}^{-1})
      \times\partial_{\xi_{n}}\Big(\frac{ic(\xi)\sum\limits_j\Big[c(dx_j)|\xi|^2
      +2\xi_{j}c(\xi)\Big]D_{x_j}(f^{-1})c(\xi)}{|\xi|^8}\Big)\Big]\nonumber\\
      &&\times(x_0)d\xi_n\sigma(\xi')dx'\nonumber\\
      &=&-\frac{15i}{2}\partial_{x_{n}}(f)\pi 8\Omega_{4}dx'.
\end{eqnarray}

Thus we have
\begin{eqnarray}
{\bf case~ (b)}
&=&(-\frac{41}{8}i-\frac{195}{8})\pi h'(0)\Omega_4dx'+120i\pi g(dx_n,\theta')\Omega_4dx'\nonumber\\
&&+\frac{3}{8f}\pi  g[dx_{n},df]\Omega_{4}dx'
      -\frac{15i}{2}\partial_{x_{n}}(f)\pi 8\Omega_{4}dx'.
\end{eqnarray}

\noindent {\bf  case (c)}~$r=-2,l=-3,|\alpha|=j=k=0$.\\

By (4.2), we have

\begin{eqnarray}
{\rm case~ (c)}&=&-i\int_{|\xi'|=1}\int^{+\infty}_{-\infty}{\rm tr} \Big[\pi^{+}_{\xi_{n}}\sigma_{-2}(f\widehat{D}_{N}^{-1})
      \times\partial_{\xi_{n}}
      \sigma_{-3}
\big(f^{-1}(\widehat{D}_{N}^{*})^{-1}\cdot f\widehat{D}_{N}^{-1}\cdot f^{-1}(\widehat{D}_{N}^{*})^{-1}\big)\Big](x_0)\nonumber\\
&&\times d\xi_n\sigma(\xi')dx'\nonumber\\
&=&-i\int_{|\xi'|=1}\int^{+\infty}_{-\infty}{\rm tr} \Big[\pi^{+}_{\xi_{n}}\sigma_{-2}(\widehat{D}_{N}^{-1})
      \times\partial_{\xi_{n}}
      \sigma_{-3}
\big((\widehat{D}_{N}^{*}\widehat{D}_{N}\widehat{D}_{N}^{*})^{-1}\big)\Big]
(x_0)d\xi_n\sigma(\xi')dx'.
\end{eqnarray}
By (4.26), we have
\begin{equation}
\partial_{\xi_{n}}\sigma_{-3}((\widehat{D}_{N}^{*}\widehat{D}_{N}\widehat{D}_{N}^{*})^{-1})=\frac{-4 i \xi_{n}c(\xi')}{(1+\xi_{n}^{2})^{3}}+\frac{i(1- 3\xi_{n}^{2})c(dx_{n})}
{(1+\xi_{n}^{2})^{3}}.
\end{equation}
By (3.46), we obtain
\begin{eqnarray}
\pi^{+}_{\xi_{n}}\Big(\sigma_{-2}(\widehat{D}_{N}^{-1})\Big)(x_{_{0}})|_{|\xi'|=1}
&=&\pi^+_{\xi_n}\Big[\frac{c(\xi)b_0^{2}(x_0)c(\xi)+c(\xi)c(dx_n)
\partial_{x_n}[c(\xi')](x_0)}{(1+\xi_n^2)^2}-h'(0)\frac{c(\xi)c(dx_n)c(\xi)}{(1+\xi_n^{2})^3}\Big]\nonumber\\
&&+\pi^+_{\xi_n}\Big[\frac{c(\xi)[b_0^{1}(x_0)]c(\xi)}{(1+\xi_n^2)^2}\Big]
+\pi^+_{\xi_n}\Big[\frac{c(\xi)[\bar{c}(\theta)+c(\theta')]c(\xi)(x_0)}{(1+\xi_n^2)^2}\Big].
\end{eqnarray}

Furthermore,
\begin{eqnarray}
&&\pi^+_{\xi_n}\Big[\frac{c(\xi)[\bar{c}(\theta)+c(\theta')](x_0)c(\xi)}
{(1+\xi_n^2)^2}\Big]\nonumber\\
&=&\pi^+_{\xi_n}\Big[\frac{c(\xi')[\bar{c}(\theta)+c(\theta')](x_0)c(\xi')}
{(1+\xi_n^2)^2}\Big]
+\pi^+_{\xi_n}\Big[ \frac{\xi_nc(\xi')[\bar{c}(\theta)+c(\theta')]
(x_0)c(dx_{n})}{(1+\xi_n^2)^2}\Big]\nonumber\\
&&+\pi^+_{\xi_n}\Big[\frac{\xi_nc(dx_{n})[\bar{c}(\theta)
+c(\theta')](x_0)c(\xi')}{(1+\xi_n^2)^2}\Big]
+\pi^+_{\xi_n}\Big[\frac{\xi_n^{2}c(dx_{n})[\bar{c}(\theta)
+c(\theta')](x_0)c(dx_{n})}{(1+\xi_n^2)^2}\Big]\nonumber\\
&=&-\frac{c(\xi')[\bar{c}(\theta)+c(\theta')](x_0)c(\xi')(2+i\xi_{n})}{4(\xi_{n}-i)^{2}}
+\frac{ic(\xi')[\bar{c}(\theta)+c(\theta')](x_0)c(dx_{n})}{4(\xi_{n}-i)^{2}}\nonumber\\
&&+\frac{ic(dx_{n})[\bar{c}(\theta)+c(\theta')](x_0)c(\xi')}{4(\xi_{n}-i)^{2}}
+\frac{-i\xi_{n}c(dx_{n})[\bar{c}(\theta)+c(\theta')](x_0)c(dx_{n})}{4(\xi_{n}-i)^{2}}.
\end{eqnarray}
By (3.47)-(3.49) and (4.54),
 we have
\begin{eqnarray}
{\rm tr }\bigg[\pi^+_{\xi_n}\Big(\frac{c(\xi)b_0^{1}(x_0)c(\xi)}{(1+\xi_n^2)^2}\Big)\times
\partial_{\xi_n}\sigma_{-3}((\widehat{D}^{*}\widehat{D}\widehat{D}^{*})^{-1})(x_0)\bigg]\bigg|_{|\xi'|=1}=\frac{2-8i\xi_n-6\xi_n^2}{4(\xi_n-i)^{2}(1+\xi_n^2)^{3}}{\rm tr }[b_0^{1}(x_0)c(\xi')],
\end{eqnarray}
By (3.52)-(3.54), we have
\begin{eqnarray}
\pi^+_{\xi_n}\Big[\frac{c(\xi)b_0^{2}(x_0)c(\xi)+c(\xi)c(dx_n)\partial_{x_n}(c(\xi'))(x_0)}{(1+\xi_n^2)^2}\Big]-h'(0)\pi^+_{\xi_n}\Big[\frac{c(\xi)c(dx_n)c(\xi)}{(1+\xi_n)^3}\Big]:= B_1-B_2,
\end{eqnarray}
where
\begin{eqnarray}
B_1&=&\frac{-1}{4(\xi_n-i)^2}\big[(2+i\xi_n)c(\xi')b^2_0c(\xi')+i\xi_nc(dx_n)b^2_0c(dx_n) \nonumber\\
&&+(2+i\xi_n)c(\xi')c(dx_n)\partial_{x_n}c(\xi')+ic(dx_n)b^2_0c(\xi')
+ic(\xi')b^2_0c(dx_n)-i\partial_{x_n}c(\xi')\big]\nonumber\\
&=&\frac{1}{4(\xi_n-i)^2}\Big[\frac{5}{2}h'(0)c(dx_n)-\frac{5i}{2}h'(0)c(\xi')
  -(2+i\xi_n)c(\xi')c(dx_n)\partial_{\xi_n}c(\xi')+i\partial_{\xi_n}c(\xi')\Big]  ;         \\
B_2&=&\frac{h'(0)}{2}\Big[\frac{c(dx_n)}{4i(\xi_n-i)}+\frac{c(dx_n)-ic(\xi')}{8(\xi_n-i)^2}
+\frac{3\xi_n-7i}{8(\xi_n-i)^3}\big(ic(\xi')-c(dx_n)\big)\Big].
\end{eqnarray}
By (4.54) and (4.60), we have
\begin{eqnarray}
&&{\rm tr }[B_2\times\partial_{\xi_n}\sigma_{-3}((\widehat{D}^{*}_{N}\widehat{D}_{N}\widehat{D}^{*}_{N})^{-1})(x_0)]|_{|\xi'|=1}\nonumber\\
&=&{\rm tr }\Big\{ \frac{h'(0)}{2}\Big[\frac{c(dx_n)}{4i(\xi_n-i)}+\frac{c(dx_n)-ic(\xi')}{8(\xi_n-i)^2}
+\frac{3\xi_n-7i}{8(\xi_n-i)^3}[ic(\xi')-c(dx_n)]\Big] \nonumber\\
&&\times\frac{-4i\xi_nc(\xi')+(i-3i\xi_n^{2})c(dx_n)}{(1+\xi_n^{2})^3}\Big\} \nonumber\\
&=&8h'(0)\frac{4i-11\xi_n-6i\xi_n^{2}+3\xi_n^{3}}{(\xi_n-i)^5(\xi_n+i)^3}.
\end{eqnarray}
Similarly, we have
\begin{eqnarray}
&&{\rm tr }[B_1\times\partial_{\xi_n}\sigma_{-3}((\widehat{D}^{*}_{N}\widehat{D}_{N}\widehat{D}^{*}_{N})^{-1})(x_0)]|_{|\xi'|=1}\nonumber\\
&=&{\rm tr }\Big\{ \frac{1}{4(\xi_n-i)^2}\Big[\frac{5}{2}h'(0)c(dx_n)-\frac{5i}{2}h'(0)c(\xi')
  -(2+i\xi_n)c(\xi')c(dx_n)\partial_{\xi_n}c(\xi')+i\partial_{\xi_n}c(\xi')\Big]\nonumber\\
&&\times \frac{-4i\xi_nc(\xi')+(i-3i\xi_n^{2})c(dx_n)}{(1+\xi_n^{2})^3}\Big\} \nonumber\\
&=&8h'(0)\frac{3+12i\xi_n+3\xi_n^{2}}{(\xi_n-i)^4(\xi_n+i)^3};\\
&&{\rm tr }\bigg[\pi^+_{\xi_n}\Big(\frac{c(\xi)[\bar{c}(\theta)+c(\theta')]
(x_0)c(\xi)}{(1+\xi_n^2)^2}\Big)\times
\partial_{\xi_n}\sigma_{-3}((\widehat{D}^{*}_{N}\widehat{D}_{N}\widehat{D}_{N}^{*})^{-1})
(x_0)\bigg]\bigg|_{|\xi'|=1}\nonumber\\
&=&\frac{2-8i\xi_n-6\xi_n^2}{4(\xi_n-i)^{2}(1+\xi_n^2)^{3}}{\rm tr }[[\bar{c}(\theta)+c(\theta')](x_0)c(\xi')]\nonumber\\
&=&\frac{2-8i\xi_n-6\xi_n^2}{4(\xi_n-i)^{2}(1+\xi_n^2)^{3}}{\rm tr }[c(\theta')(x_0)c(\xi')]\nonumber\\
&=&\frac{2-8i\xi_n-6\xi_n^2}{4(\xi_n-i)^{2}(1+\xi_n^2)^{3}}[-g(\theta',\xi')]
{\rm tr}[{\rm id}].
\end{eqnarray}

By $\int_{|\xi'|=1}\xi_{1}\cdot\cdot\cdot\xi_{2q+1}\sigma(\xi')=0,$ we have
\begin{eqnarray}
{\rm case~(c)}&=&
 -i h'(0)\int_{|\xi'|=1}\int^{+\infty}_{-\infty}
 8\times\frac{-7i+26\xi_n+15i\xi_n^{2}}{(\xi_n-i)^5(\xi_n+i)^3}d\xi_n\sigma(\xi')dx' \nonumber\\
 &&-i\int_{|\xi'|=1}\int^{+\infty}_{-\infty}
 \bigg[\frac{2-8i\xi_n-6\xi_n^2}{4(\xi_n-i)^{2}(1+\xi_n^2)^{3}}[-g(\theta',\xi')]
{\rm tr}[{\rm id}]\bigg]\bigg|_{|\xi'|=1}d\xi_n\sigma(\xi')dx'\nonumber\\
&=&-8i h'(0)\times\frac{2 \pi i}{4!}\Big[\frac{-7i+26\xi_n+15i\xi_n^{2}}{(\xi_n+i)^3}
     \Big]^{(5)}|_{\xi_n=i}\Omega_4dx'\nonumber\\
&=&\frac{55}{2}\pi h'(0)\Omega_4dx'.
\end{eqnarray}

Now $\Psi$ is the sum of the {\bf  case (a)}, {\bf  case (b)} and {\bf  case (c)}, then
\begin{eqnarray}
\Psi&=&(10i+88)\pi f^{-1}\partial_{x_{n}}(f)\cdot{\rm \Omega_{4}}dx'+\frac{\pi i}{2}\cdot f\cdot\partial_{x_{n}}(f^{-1})\Omega_{4}dx'-(\frac{41}{8}i-\frac{65}{8})\pi h'(0)\Omega_4dx'\nonumber\\
&&+120i\pi g(dx_n,\theta')\Omega_4dx'+\frac{3}{8f}\pi  g[dx_{n},df]\Omega_{4}dx'
      -60i\partial_{x_{n}}(f)\pi \Omega_{4}dx'.
\end{eqnarray}
\begin{thm}
Let $M$ be $6$-dimensional oriented
compact manifolds with the boundary $\partial M$ and the metric
$g^M$ as above, $\widehat{D}_{N}$ and $\widehat{D}^*_{N}$ be modified Novikov operators on $\widehat{M}$, then
 \begin{eqnarray}
&&\widetilde{Wres}\big[\pi^{+}(f\widehat{D}_{N}^{-1}) \circ\pi^{+}\big(f^{-1}(\widehat{D}_{N}^{*})^{-1}\cdot f\widehat{D}_{N}^{-1}\cdot f^{-1}(\widehat{D}_{N}^{*})^{-1}\big)\big]\nonumber\\
&=&128\pi^{3}\int_{M}
2^6\bigg\{\bigg(-\frac{1}{12}s-|\theta|^2+(n-2)|\theta'|^2+g(\widetilde{e_{j}},\nabla^{TM}_{\widetilde{e_{j}}}\theta')\bigg)+\bigg[f^{-1}\Delta(f)
+\langle grad_{M}f,grad_{M}f^{-1}\rangle\bigg]\nonumber\\
&&-\Big(\bar{c}(\theta)-c(\theta')\Big)c(df)f^{-1}
-\frac{1}{2}\bigg[f^{-1}\Delta(f)(x_0)+\langle grad_{M}f,grad_{M}f^{-1}\rangle\bigg]+\frac{1}{4}\sum_{i}\Big[c(e_{i})c(\theta')-c(\theta')
c(e_{i})\Big]\nonumber\\
&&\times c(e_{i})c(df)f^{-1}+\frac{1}{4}\sum_{i}c(e_{i})c(df)f^{-1}
\Big[c(e_{i})c(\theta')-c(\theta')c(e_{i})\Big]\bigg\}d{\rm Vol_{M}}+
       \int_{\partial M}\Bigg\{(10i+88)\pi f^{-1}\partial_{x_{n}}(f)\nonumber\\
&&\times{\rm \Omega_{4}}dx'+\frac{\pi i}{2}\cdot f\cdot\partial_{x_{n}}(f^{-1})\Omega_{4}dx'(-\frac{41}{8}i+\frac{65}{8})\pi h'(0)\Omega_4dx'+120i\pi g(dx_n,\theta')\Omega_4dx'+\frac{3}{8f}\pi  g(dx_{n},df)\Omega_{4}dx'\nonumber\\
&&-60i\partial_{x_{n}}(f)\pi \Omega_{4}dx'
      \Bigg\}dvol_{M},
\end{eqnarray}
where  $s$ is the scalar curvature.
\end{thm}

\section*{Acknowledgements}
The first author is supported by DUFE202159. The
partial research of the corresponding author was supported by NSFC. 11771070.

\end{document}